\documentclass[twoside,fleqn]{article}
\usepackage{amssymb,amsfonts,amsthm}
\makeatletter
\def\fileversion{v2.6}
\def\filedate{24 November 1993}

\newdimen\@bls                    % \@b(ase)l(ine)s(kip)
\@bls=\baselineskip               % \@bls ~= \baselineskip for \normalsize
\advance\@bls -1ex                % (fudge term)
\newdimen\@eps                    %
\def\section{\@startsection{section}{1}{\z@}
    {1.5\@bls plus 0.5\@bls}{1\@bls}{\normalsize\bf}}
\def\subsection{\@startsection{subsection}{2}{\z@}
    {1\@bls plus 0.25\@bls}{\@eps}{\normalsize\bf}}
\def\subsubsection{\@startsection{subsubsection}{3}{\z@}
    {1\@bls plus 0.25\@bls}{\@eps}{\normalsize\bf}}
\def\paragraph{\@startsection{paragraph}{4}{\parindent}
    {1\@bls plus 0.25\@bls}{0.5em}{\normalsize\bf}}
\def\subparagraph{\@startsection{subparagraph}{4}{\parindent}
    {1\@bls plus 0.25\@bls}{0.5em}{\normalsize\bf}}

\def\@sect#1#2#3#4#5#6[#7]#8{\ifnum #2>\c@secnumdepth
    \def\@svsec{}\else
    \refstepcounter{#1}\edef\@svsec{\csname the#1\endcsname.\hskip0.5em}\fi
    \@tempskipa #5\relax
    \ifdim \@tempskipa>\z@
      \begingroup
        #6\relax
        \@hangfrom{\hskip #3\relax\@svsec}{\interlinepenalty \@M #8\par}%
      \endgroup
      \csname #1mark\endcsname{#7}\addcontentsline
        {toc}{#1}{\ifnum #2>\c@secnumdepth \else
          \protect\numberline{\csname the#1\endcsname}\fi #7}%
    \else
      \def\@svsechd{#6\hskip #3\@svsec #8\csname #1mark\endcsname
        {#7}\addcontentsline{toc}{#1}{\ifnum #2>\c@secnumdepth \else
          \protect\numberline{\csname the#1\endcsname}\fi #7}}%
    \fi \@xsect{#5}}

\long\def\@makefigurecaption#1#2{\vskip 10mm #1. #2\par}

\long\def\@maketablecaption#1#2{\hbox to \hsize{\parbox[t]{\hsize}
    {#1 \\ #2}}\vskip 0.3ex}

\def\fnum@figure{Figure \thefigure}
\def\figure{\let\@makecaption\@makefigurecaption \@float{figure}}
\@namedef{figure*}{\let\@makecaption\@makefigurecaption \@dblfloat{figure}}

\def\table{\let\@makecaption\@maketablecaption \@float{table}}
\@namedef{table*}{\let\@makecaption\@maketablecaption \@dblfloat{table}}

\textfloatsep=\floatsep           % Space between main text and floats
\intextsep=\floatsep              % Space between in-text figures and
\long\def\@makefntext#1{\parindent 1em\noindent\hbox{${}^{\@thefnmark}$}#1}

\def\maketitle{\begingroup        % Initialize generation of front-matter
      \def\thefootnote{\fnsymbol{footnote}}%
      \newpage \global\@topnum\z@
      \@maketitle \@thanks
    \endgroup
    \let\maketitle\relax \let\@maketitle\relax
    \gdef\@thanks{}\let\thanks\relax
    \gdef\@address{}\gdef\@author{}\gdef\@title{}\let\address\relax}

\def\justify@on{\let\\=\@normalcr
    \leftskip\z@ \@rightskip\z@ \rightskip\@rightskip}

\newbox\fm@box                    % Box to capture front-matter in

\def\@maketitle{%                 % Actual formatting of \maketitle
    \global\setbox\fm@box=\vbox\bgroup
      \vskip 8mm                    % 930715: 8mm white space above title
      \raggedright                  % Front-matter text is ragged right
      \hyphenpenalty\@M             % and is not hyphenated.
      {\Large \@title \par}         % Title set in larger font.
      \vskip\@bls                   % One line of vertical space after title.
      {\normalsize                  % each author set in the normal
       \@author \par}               % typeface size
      \vskip\@bls                   % One line of vertical space after
author(s).
      \@address                     % all addresses
    \egroup
    \twocolumn[%                    % Front-matter text is over 2 columns.
      \unvbox\fm@box                % Unwrap contents of front-matter box
      \vskip\@bls                   % add 1 line of vertical space,
      \unvbox\abstract@box          % unwrap contents of abstract boxes,
      \vskip 2pc]}                  % and add 2pc of vertical space

\newcounter{address}
\def\theaddress{\alph{address}}
\def\@makeadmark#1{\hbox{$^{\rm #1}$}}

\def\address#1{\addressmark\begingroup
    \xdef\@tempa{\theaddress}\let\\=\relax
    \def\protect{\noexpand\protect\noexpand}\xdef\@address{\@address
    \protect\addresstext{\@tempa}{#1}}\endgroup}
\def\@address{}

\def\addressmark{\stepcounter{address}%
    \xdef\@tempb{\theaddress}\@makeadmark{\@tempb}}

\def\addresstext#1#2{\leavevmode \begingroup
    \raggedright \hyphenpenalty\@M \@makeadmark{#1}#2\par \endgroup
    \vskip\@bls}

\newbox\abstract@box              % Box to capture abstract in

\def\abstract{%
    \global\setbox\abstract@box=\vbox\bgroup
    \small\rm
    \ignorespaces}
\def\endabstract{\par \egroup}

\def\thebibliography#1{\section*{REFERENCES}\list{\arabic{enumi}.}
    {\settowidth\labelwidth{#1.}\leftmargin=1.67em
     \labelsep\leftmargin \advance\labelsep-\labelwidth
     \itemsep\z@ \parsep\z@
     \usecounter{enumi}}\def\makelabel##1{\rlap{##1}\hss}%
     \def\newblock{\hskip 0.11em plus 0.33em minus -0.07em}
     \sloppy \clubpenalty=4000 \widowpenalty=4000 \sfcode`\.=1000\relax}

\newcount\@tempcntc
\def\@citex[#1]#2{\if@filesw\immediate\write\@auxout{\string\citation{#2}}\fi
    \@tempcnta\z@\@tempcntb\m@ne\def\@citea{}\@cite{\@for\@citeb:=#2\do
      {\@ifundefined
         {b@\@citeb}{\@citeo\@tempcntb\m@ne\@citea
          \def\@citea{,\penalty\@m\ }{\bf ?}\@warning
         {Citation `\@citeb' on page \thepage \space undefined}}%
      {\setbox\z@\hbox{\global\@tempcntc0\csname b@\@citeb\endcsname\relax}%
       \ifnum\@tempcntc=\z@ \@citeo\@tempcntb\m@ne
         \@citea\def\@citea{,\penalty\@m}
         \hbox{\csname b@\@citeb\endcsname}%
       \else
        \advance\@tempcntb\@ne
        \ifnum\@tempcntb=\@tempcntc
        \else\advance\@tempcntb\m@ne\@citeo
        \@tempcnta\@tempcntc\@tempcntb\@tempcntc\fi\fi}}\@citeo}{#1}}

\def\@citeo{\ifnum\@tempcnta>\@tempcntb\else\@citea
    \def\@citea{,\penalty\@m}%
    \ifnum\@tempcnta=\@tempcntb\the\@tempcnta\else
     {\advance\@tempcnta\@ne\ifnum\@tempcnta=\@tempcntb \else
\def\@citea{--}\fi
      \advance\@tempcnta\m@ne\the\@tempcnta\@citea\the\@tempcntb}\fi\fi}

\def\ps@crcplain{\let\@mkboth\@gobbletwo
       \def\@oddhead{\reset@font{\sl\rightmark}\hfil \rm\thepage}%
       \def\@evenhead{\reset@font\rm \thepage\hfil\sl\leftmark}%
       \let\@oddfoot\@empty
       \let\@evenfoot\@oddfoot}

\ps@crcplain                    % modified 'plain' page style

\theoremstyle{plain} %% This is the default
\newtheorem{thm}{Theorem}[section]
\newtheorem{cor}[thm]{Corollary}
\newtheorem{lem}[thm]{Lemma}
\newtheorem{prop}[thm]{Proposition}

\theoremstyle{definition}

\newtheorem{rem}[thm]{Remark}

\theoremstyle{remark}

\def\noi{\noindent}

\def\tensor{\otimes}
\def\too{\longrightarrow}

\def\={\cong}
\def\>{\supset}
\def\<{\subset}

\def\12{\frac{1}{2}}
\def\0{^{\circ}}

\def\CC{{\bf C}}

\def\NN{{\bf N}}

\def\RR{{\bf R}}

\def\ZZ{{\bf Z}}

\def\Aa{{\mathcal A}}
\def\Bb{{\mathcal B}}
\def\Cc{{\mathcal C}}
\def\Dd{{\mathcal D}}
\def\Ee{{\mathcal E}}
\def\Ff{{\mathcal F}}

\def\Hh{{\mathcal H}}
\def\Kk{{\mathcal K}}
\def\Ll{{\mathcal L}}

\def\Nn{{\mathcal N}}

\def\Pp{{\mathcal P}}

\def\Rr{{\mathcal R}}
\def\Ss{{\mathcal S}}

\def\={\cong}
\def\>{\supset}
\def\<{\subset}

\def\12{\frac{1}{2}}
\def\2{\Dd}
\def\3{\Nn}
\def\4{\Rr}
\def\6{\cup}
\def\8{\otimes}
\def\0{^{\circ}}

\def\a{\alpha}

\def\C{\CC}
\def\e{\varepsilon}
\def\g{\gamma}
\def\G{\Gamma}

\def\k{\kappa}

\def\la{\lambda}

\def\m{\mu}

\def\N{\NN}
\def\p{\pi}

\def\R{\RR}
\def\s{\sigma}
\def\Si{\Sigma}
\def\z{\zeta}
\def\Z{\ZZ}

\def\Si{S\kern -.65em /}

\def\tensor{\otimes}

\def\Tr{\mbox{\rm Tr\,}}

\makeatother

\newcommand{\bee}{\begin{equation}}
\newcommand{\ene}{\end{equation}}

\title%[Analytic Surgery of the $\z$-determinant]
{Analytic Surgery of the $\z$-determinant of the Dirac operator}

\author{{Jinsung Park, Krzysztof P. Wojciechowski}
\address{Department of Mathematics\\IUPUI (Indiana/Purdue)\\
Indianapolis IN 46202--3216, U.S.A.}
\thanks{First author partially supported by Korea Science and
Engineering Foundation}}

\begin{document}

\begin{abstract}
We review the work of the authors and their collaborators on the decomposition
of the $\z$-determinant of the Dirac operator into the contributions
coming from different parts of a manifold.
\end{abstract}

\maketitle

\section{Introduction}\label{s:intr}

\bigskip

The main theme of our lectures is to discuss how the
decomposition of a manifold (space-time) affects the structure
of the $\z$-determinant, which is a delicate spectral invariant.
This subject has been studied by many authors from many different
perspectives (see for instance \cite{BFK92}, \cite{CLM1},
\cite{Ch83}, \cite{Ch87}, \cite{HMaM95}, \cite{MaM95}, \cite{Mu94}
, \cite{Mu98}, \cite{P96}, \cite{P98} and infinitely many others).
They have used many different technical approaches introducing
incredible amount of beautiful and difficult mathematics. These
notes are meant to be an introduction to the authors' perspective
onto the subject. The focus here is on ideas rather than on
rigorous arguments. Most of the results have been published in
recent papers by the authors and their collaborators and we give
precise bibliographical references. However, let us stress that
due to enormously rich literature we do not attempt to be as
complete as possible. We want to apologize for not mentioning many
important works, that have made an enormous impact on this area of
mathematics and mathematical physics.

\bigskip

In Section 2 we study the properties of the $\z$-determinant of the Dirac
operator on a closed manifold
using  the Heat Equation method. We present here standard material, that is
described in many great sources.
In Section 3 we describe the adjustment we have to make in order
to study Dirac  operators on a manifold with boundary. We explain our
choice of the space of the boundary
conditions and show that there is a natural notion of the determinant
related to this space.
We discuss the projective equality of this new determinant
to the $\z$-determinant of the boundary problems for the Dirac operators
established in the recent work of Scott and Wojciechowski
(see \cite{SSKPW299}, see also \cite{SSKPW00} for the additional discussion).
In Section 4 we outline our method of analyzing the decomposition of
the $\z$-determinant.
Section 5 deals with the boundary contributions which appear when we split
a manifold along the submanifold of codimension $1$.
Then in Section 6 we explain how to use the adiabatic approach in order to
separate the contributions coming from different
parts of the manifold and the boundary contributions. This decomposition is
completed in Section 7.
In Section 8 we present the ``adiabatic'' decomposition formulas for the
$\z$-determinant of the Dirac Laplacians.
Let us point out that formulas (\ref{e:d3}) and (\ref{e:n1})  are new,
while the complete proof of the formula (\ref{e:chsplit1}) was
given  in a recent paper by the authors (see \cite{JKPW2}, see also
\cite{JPKW1}).
In Section 9 we discuss the decomposition of the ``phase"of the
$\z$-determinant, the $\eta$-invariant.
Here we make more comments concerning the analysis on a manifold with
boundary. We explain
why there are no analytical problems with the definition of the
$\z$-determinant on $Gr_{\infty}^*(\Dd)$,
the Grassmannian of the boundary conditions we discuss in this paper. Then
we present the proof of the decomposition
formula for the $\eta$-invariant. The disadvantage of our method is that it
does not tell us anything
about the integer contribution. Additional study is needed to
detect the integer contribution, which is responsible for
some intriguing topological phenomena. Due to the
lack of expertise and space in this article we
do not discuss this topic. Instead of that, we refer to a beautiful, recent work
of Kirk and Lesch \cite{KL00}.
In the last Section, we discuss the invariance of the ratio  of the
$\z$-determinants of two elliptic problems
with respect to the length of the collar neighborhood of the boundary. It
is well-known, that in general, the $\z$-determinant
changes when we stretch the collar. We discuss here the case in which the
ratio of the determinants of
two Atiyah--Patodi--Singer problems remains constant. The proof is based on
the results of the work of Scott and Wojciechowski discussed
in Section 3 (see \cite{SSKPW299}).

\bigskip

In the reminder of the Introduction we introduce the main hero of the
lectures -
the $\z$-determinant of the Dirac operator. We follow here a beautiful
exposition
given by Singer in \cite{Si85}.

In many important problems of quantizing gauge theories, as well
as in some mathematical problems, it is necessary to discuss
directly a regularized determinant of an elliptic operator. The
{\it Heuristic Approach} to the determinant in this context was
first proposed by mathematicians for the case of a positive
definite second-order elliptic differential operator

$$L: C^{\infty}(M;S) \to C^{\infty}(M;S)$$

\noi
acting on sections of a smooth vector bundle $S$ over a
closed manifold $M$. The operator $L$ has a discrete spectral
resolution and therefore formally has determinant equal to the infinite
product of its eigenvalues. The starting point in defining a
regularized product is the following formula for an invertible
finite-rank linear operator T:

\begin{equation}\label{e:03}
{\rm ln} \; {\rm det} \ T = - {d\over{ds}}\{ {\rm Tr} \
T^{-s}\}|_{s=0} \ \ .
\end{equation}

\noi
For large $Re(s)$ the $\z$-function of the operator $L$ is
just the trace occurring on the right side of (\ref{e:03})

\begin{equation}\label{e:04}
\z_{L}(s) = {\rm Tr} \ L^{-s} =
{1\over{\G(s)}}\int_{0}^{\infty}t^{s-1}{\rm Tr} \ e^{-tL}dt .
\end{equation}

\bigskip

\noi It is a holomorphic function of $s$ for $Re(s) > {{dim \
M}\over{2}}$ and has a meromorphic extension to the whole complex
plane with only simple poles (see \cite{Se66}). In particular
$s=0$ is not a pole. Hence $\z_{L}'(0) =
{d\over{ds}}\{\z_{L}(s)\}|_{s=0}$ is well-defined and we may
define the $\z$-{\em determinant} by

\begin{equation}\label{e:05}
{\rm det}_{\z}L = e^{-\z_{L}'(0)} \ \, .
\end{equation}

\noi This definition was introduced in 1971, in  a famous paper of
Ray and Singer \cite{RaSi71},  in order to define {\it Analytic
Torsion}, the analytical counterpart to the topological invariant
{\it Franz-Reidemeister Torsion}. The equality of the two torsions
was subsequently proved independently by Jeff Cheeger and Werner
M$\ddot {\rm u}$ller (see \cite{Ch79}, \cite{Mu78}). Since then,
there have been numerous applications of the $\z$-determinant in
physics and mathematics, beginning with the 1977 Hawking paper
\cite{Ha75} on quantum gravity.

For positive-definite operators of Laplace type over a closed
manifold the $\z$-determinant provides a generally satisfactory
regularization method. Though the fundamental multiplicative
property of the determinant no longer holds; if $L_1$ and $L_2$
denote two positive elliptic operator of
positive order on a Hilbert space $H$ then in general

$$
{\rm det}_{\z}L_1L_2 \ne {\rm det}_{\z}L_1{\cdot}
{\rm det}_{\z}L_2 \ \ .$$

\bigskip

\noi
We refer to other talks in the Meeting for a discussion of the
so-called Multiplicative Anomaly. In many physical applications,
however, such as the quantization of Fermions, one encounters the
more problematic task of defining the determinant of a first-order
{\it Dirac operator}. These are not positive operators, and now
the gauge anomalies may arise due to the phase of
the determinant (see \cite{AtSi84}). For a Dirac operator
$\Dd: C^{\infty}(M;S) \to C^{\infty}(M;S)$ acting on sections of
a bundle of Clifford modules over a closed (odd-dimensional)
manifold $M$ one proceeds in the way outlined below. The operator $\Dd$ is an
elliptic self-adjoint first-order operator and hence has
infinitely many positive and negative eigenvalues. Let
$\{\la_k\}_{k \in \N}$ denote the set of positive eigenvalues and
$\{-\mu_k\}_{k \in \N}$ denote the set of negative eigenvalues.
Once again, $\z_{\Dd}(s) = \Tr(\Dd^{-s})$ is well-defined and
holomorphic for ${\rm Re}(s) > {\rm dim} \ M$ and we have

$$\z_\Dd(s) = \sum_k \la_k^{-s} + \sum_k (-1)^{-s}\mu_k^{-s} $$
$$=\sum_k({{\la_k^{-s} - \mu_k^{-s}}\over{2}} + {{\la_k^{-s} +
\mu_k^{-s}}\over{2}})$$
$$\qquad + (-1)^{-s}\sum_k({{\la_k^{-s} +
\mu_k^{-s}}\over{2}} - {{\la_k^{-s} - \mu_k^{-s}}\over{2}}) \ \
,$$

\bigskip

\noi
which can be written as

\begin{equation}\label{e:06}
\z_\Dd(s) = (-1)^{-s}{{\z_{\Dd^2}(s/2) - \eta_{\Dd}(s)}\over{2}}
\ \,
\end{equation}

$$+{{ \z_{\Dd^2}(s/2)+\eta_{\Dd}(s)}\over{2}} \ , \qquad$$

\bigskip

\noi
where $\eta_{\Dd}(s) = \sum_k{\la_k^{-s}} -
\sum_k{\mu_k^{-s}}$ is the $\eta$-function of the operator $\Dd$
introduced by Atiyah, Patodi and Singer (see \cite{AtPaSi75}).
Once again, it is holomorphic for $Re(s)$ large and has a
meromorphic extension to the whole complex plane with only simple
poles. There is no pole at $s=0$ and therefore we can study the
derivative of $\z_{\Dd}(s)$ at $s=0$. We have

$$\z'_{\Dd}(0) = {{\z'_{\Dd^2}(0)}\over{2}} +
{d\over{ds}}\{(-1)^{-s}\}|_{s=0} {{\z_{\Dd^2}(0) -
\eta_{\Dd}(0)}\over{2}} \ .$$

\bigskip

\noi
The ambiguity in defining $(-1)^{-s}$ (i.e. a choice of
spectral cut) now leads to an ambiguity in the phase of the
$\z$-determinant. We have

$$(-1)^{-s} = e^{{\pm}i{\p}s} \ \ ,$$

\noi
and we pick the $``-"$ sign. This leads to the following
formula for the $\z$-determinant of the Dirac operator $\Dd$:

\begin{equation}\label{e:07}
{\rm det}_{\z}\Dd = e^{{{i\p}\over{2}}(\z_{\Dd^2}(0) -
\eta_{\Dd}(0))}{\cdot} e^{-{1\over{2}}\z'_{\Dd^2}(0)} \ \, .
\end{equation}

\bigskip

\begin{rem}
We refer to Section 7 of \cite{SSKPW299} for a discussion of the
choice of sign of the phase of the $\z$-determinant.
\end{rem}

\bigskip
We need to study more closely the regularization process used to make the
definition (\ref{e:07}). This will be done in the next Section,
where the Heat Equation enters the scene.

\bigskip

\section{$\z$-determinant and Heat Equation}\label{s:he}

\bigskip

We use the {\it Heat Equation} method to make sense of the
$\z$-determinant. We recall the standard material (see \cite{Gi95}
for details). In this Section, we assume that $\Dd$ has the
trivial kernel for convenience. The key are the following
formulas:

\begin{equation}\label{e:h4}
\z_{\Dd^2}(s) = {\rm Tr} (D^2)^{-{s}} =
{1\over{\G(s)}}\int_0^{\infty}t^{s-1}{\rm Tr} \ e^{-t\Dd^2}dt \ \,
\end{equation}

\bigskip

\noi
for $Re(s) > {{dim \ M}\over{2}}$ and

$$\eta_{\Dd}(s) = {\rm Tr} \ D(D^2)^{-{{s-1}\over{2}}}
$$
$$ \qquad
={1\over{\G({{{s+1}\over{2}})}}}\int_0^{\infty}t^{{{s-1}\over{2}}}{\rm
Tr} \ {\Dd}e^{-t\Dd^2}dt \ ,$$

\bigskip

\noi for $Re(s) > {{1 + dim \ M}\over{2}}$.

\vskip 5mm

We prove the second equality in (\ref{e:h4}). The proof of the first one is
completely analogous. We have
$$
\int_0^{\infty}t^{{{s-1}\over{2}}}{\rm Tr} \ {\Dd}e^{-t\Dd^2}dt =
\sum_{-\infty}^{+\infty}
\int_0^{\infty}t^{{{s-1}\over{2}}}\la_ke^{-t\la_k^2}dt
$$
$$
=\sum_{-\infty}^{+\infty}\la_k(\la_k^2)^{-{{s+1}\over{2}}}\int_0^{\infty}
(t\la_k^2)^{{{s-1}\over{2}}}e^{-t\la_k^2}d(t\la_k^2)
$$
$$
=\sum_{-\infty}^{+\infty}{\rm sign} \ \la_k{\cdot}|\la_k|^{-s}{\cdot}
\int_0^{\infty}r^{{{s-1}\over{2}}}e^{-r}dr
$$
$$
=\G \left({{{s+1}\over{2}}}\right) \eta_{\Dd}(s) .$$

\bigskip
These formulas hold for $s$ making the operators $(D^2)^{-{s}}$
and $D(D^2)^{-{{s-1}\over{2}}}$ operators of trace class.
Now we expand the $\z$-function and $\eta$-function to the whole
complex plane. We use here the well-known fact that the trace ${\rm Tr} \
e^{-t\Dd^2}$ has an asymptotic expansion of the form

\begin{equation}\label{e:expansion0}
{\rm Tr} \ e^{-t\Dd^2} = t^{-{n\over{2}}}\sum_{k=0}^N t^ka_k +
O(t^{N + 1 - {\frac {n}{2}}}) \ \, .
\end{equation}

\noi A more general formula (proved in \cite{Gi95} Section 1.9.)
gives the following expansion:
\begin{equation}\label{e:expansion1}
{\rm Tr} \ A e^{-t\Dd^2} = \sum_{k=0}^N t^{{{k - n - a}\over{2}}}b_k +
O(t^{N + 1 -{\frac {a + n}{2}}}) \ \,
\end{equation}

\noi
where $A$ denotes a differential operator of order $a$. The coefficients
$a_k$ and $b_k$ are the integrals of the local densities

$$a_k = {\int_M}{\alpha}_k(x)dx \ \ {\rm and}
\ \ b_k = {\int_M}{\beta}_k(x)dx \ ,$$

\noi
where $\alpha_k(x)$ is constructed from the coefficients of $\Dd$ at the
point $x \in M$ and
$\beta_k(x)$ is constructed from coefficients of $A$ and $\Dd$ at $x$ .
Moreover,

$$
\beta_k(x) = 0 \ \ {\rm for} \ \ k + a \ \ {\rm odd} \ .$$

\bigskip

Now we see how to extend $\z_{\Dd^2}(s)$ to the whole complex plane.
$$
\int_0^{\infty}t^{s-1}{\rm Tr} \ e^{-t\Dd^2}dt
$$
$$
=\int_0^1 t^{s-1}{\rm Tr} \ e^{-t\Dd^2}dt
+ \int_1^{\infty}t^{s-1}{\rm Tr} \ e^{-t\Dd^2}dt
$$
$$
=\int_0^1t^{s-1}t^{-{n\over{2}}} \sum_{k=0}^N t^{k}a_{k}dt +
\int_1^{\infty}t^{s-1}{\rm Tr} \ e^{-t\Dd^2}dt
$$
$$
+ O(t^{s+N + 1 - {\frac {n}{2}}}) \ \ .
$$
\noi The second and the third term on the right side above provide
us with $h$, a holomorphic function of $s$ for $Re(s) >
{n\over{2}} - N -1$ and we obtain

\begin{equation}\label{e:zas}
\int_0^{\infty}t^{s-1}{\rm Tr} \ e^{-t\Dd^2}dt = \sum_{k=0}^N
{{a_k}\over{s + k -{n\over{2}}}} + h(s). \ \,
\end{equation}

\noi
It follows that $\z_{\Dd^2}(s)$ has a meromorphic extension to the whole
complex plane $\C$ with simple poles
at $s_k = {n\over{2}} - k$ , with residue equal to

$$Res_{s = {n\over{2}} - k}\z_{\Dd^2}(s) = {a_k\over{\G({n\over{2}} - k})}
\ \ .$$

\noi Let us observe a simple corollary of this analysis:

\bigskip

\begin{lem}\label{l:z}
The point $s=0$ is never a pole and $\z_{\Dd^2}(0) = 0$ for $n$ odd, and it is
equal to $a_{{n\over{2}}}$ for $n$ even.
\end{lem}

\bigskip

The reason for the regularity here is that in the neighborhood of $s=0$ ,
$\z_{\Dd^2}(s)$ can be represented in the form

\begin{equation}\label{e:zreg}
\z_{\Dd^2}(s) = {1\over{\G(s)}} \left({{a_{{n\over{2}}}}\over{s}}
+ h_1(s)\right),
\end{equation}

\noi
where $h_1$ is holomorphic in a neighborhood of $s=0$ , and the
singularity vanishes since

$$\G(s) = {1\over{s}} + \g + s{\cdot}h_2(s) \ \ ,$$

\noi
where $h_{2}$ is a holomorphic function near $s=0$ and $\g$ denotes the
Euler constant. Unfortunately this is not the case when we discuss
the $\eta$-function. The pole of ${\rm Tr} \ \Dd e^{-t\Dd^2}$ is
not cancelled out by the corresponding
pole of $\G(s)$. A more subtle argument has to be used. However the result
holds and in fact is true.

\bigskip

\begin{thm}\label{t:eta}
(see \cite{BiFr86} and \cite{Gi95}). Let $\eta_{\Dd}(s;x)$ denote the
local $\eta$-density

$$\eta_{\Dd}(s;x) = {1\over{\G({{{s+1}\over{2}})}}}\int_0^{\infty}t^{s-1}
{\rm tr}
\ \Ff(t;x,x)dt \ \ ,$$

\noi
where $\Ff(t;x,y)$ denotes the kernel of the operator $\Dd e^{-t\Dd^2}$ .
For each $x \in M$ the function
$\eta_{\Dd}(s;x)$ is a holomorphic function of $s$ for $Re(s) > -2$ .
\end{thm}

\bigskip

\begin{rem}\label{r:eta}

(1) One can view this result as the  odd-dimensional variant of the {\it
``Local Index Theorem"} for
{\it compatible Dirac operators}.

(2) It follows that the following equality holds for any compatible Dirac
operator:

\begin{equation}\label{e:ata}
\eta_{\Dd}(0) =
{1\over{\sqrt{\pi}}}\int_0^{\infty}{1\over{\sqrt{t}}}{\rm Tr} \ \Dd
e^{-t\Dd^2}dt. \ \,
\end{equation}
\end{rem}

\bigskip

To get a useful local invariant out of the $\eta$-function , we
have to study the variation  of the $\eta$-invariant (i.e.
$\frac{d}{dr}\eta_{\Dd_r}(0)$). Let us assume that
$\{\Dd_r\}_{(-\e,+\e)}$ is a smooth family of compatible Dirac
operators. For simplicity we also assume that $\Dd_r$ is an
invertible operator for any $r$ . We have to differentiate the
{\it Heat Operator} $e^{-t\Dd^2}$ . In order to do this we
introduce {\it Duhamel's Principle}.

\bigskip

{\bf Duhamel's Principle}

\noi Let $A$ and $B$ denote self-adjoint operators acting on a
separable Hilbert space $\Hh$ . The following equality holds
(under the appropriate technical assumptions):

\begin{equation}\label{e:dh1}
e^{-tA} - e^{-tB} = \int_0^te^{-sA}(B-A)e^{-(t-s)B}ds. \ \,
\end{equation}

\bigskip
An immediate consequence of Duhamel's principle that
we need is the following Proposition:
\bigskip

\begin{prop}\label{p:dh1}
The following equality holds:

\begin{equation}\label{e:dh2}
\ \ \ \ \ \ \ \ \ \ \ \ \ \ \
 {d\over{dr}}\{{\rm Tr} \ \Dd_r
e^{-t\Dd_r^2}\}|_{r=0} \ \,
\end{equation}
$$  ={\rm Tr} \ {\dot \Dd}_0e^{-t\Dd_0^2}- 2t{\cdot} \ {\rm Tr} \ {\dot
\Dd}_0\Dd_0^2e^{-t\Dd_0^2} \ \ ,$$ \noi where ${\dot \Dd}_0 =
{\frac d{dr}}\Dd_r|_{r=0}$ .
\end{prop}

\bigskip

\begin{proof}
We have

$$
{d\over{dr}}\{{\rm Tr} \ \Dd_r e^{-t\Dd_r^2}\}|_{r=0} =
\lim_{\delta \to 0} {\rm Tr}{{\Dd_{\delta} e^{-t\Dd_{\delta}^2} -
\Dd_0 e^{-t\Dd_0^2}}\over{\delta}}
$$
$$
= \lim_{\delta \to 0} {\rm Tr}{{\Dd_{\delta} -
\Dd_0}\over{\delta}}e^{-t\Dd_{\delta}^2} + \lim_{\delta \to 0}{\rm
Tr} \ \Dd_0{{e^{-t\Dd_{\delta}^2} - e^{-t\Dd_0^2}}\over{\delta}}
$$
$$
= {\rm Tr} \ {\dot
\Dd}_0e^{-t\Dd_0^2}\qquad\qquad\qquad\qquad\qquad\qquad\qquad\qquad
$$
$$
+\lim_{\delta \to 0}\{{\rm Tr} \ \Dd_0
\int_0^te^{-s\Dd_0^2}{{\Dd_0^2 -
\Dd_{\delta}^2}\over{\delta}}e^{-(t-s)\Dd_0^2}ds
$$
$$+ {\rm Tr} \ \Dd_0 \int_0^t{{(e^{-s\Dd_{\delta}^2} -
e^{-s\Dd_0^2})(\Dd_0^2 -
\Dd_{\delta}^2)}\over{\delta}}e^{-(t-s)\Dd_0^2}ds\} \ .
$$

\bigskip

\noi
The last term on the right side is of order $O(\delta)$ and we obtain

$$
{\rm Tr} \ {\dot \Dd}_0e^{-t\Dd_0^2}\qquad\qquad\qquad\qquad\qquad
$$
$$ - {\rm Tr} \ \Dd_0 \int_0^te^{-s\Dd_0^2}({\dot \Dd}_0\Dd_0 +\Dd_0{\dot
\Dd}_0)e^{-(t-s)\Dd_0^2}ds$$
$$= {\rm Tr} \ {\dot \Dd}_0e^{-t\Dd_0^2}
- 2t{\cdot}{\rm Tr}  \ {\dot \Dd}_0\Dd_0^2e^{-t\Dd_0^2} \ .
$$

\end{proof}

\bigskip
We can now discuss two formulas for the variation of the $\eta$-invariant.
The first follows from formulas (\ref{e:ata}) and (\ref{e:dh2}). We have

$${d\over{dr}}\{\eta_{\Dd_r}(0)\}|_{r=0} = {1\over{\sqrt{\pi}}}
\left \{\int_0^{\infty}{1\over{\sqrt{t}}}{\rm Tr} \ {\dot \Dd}_0e^{-t\Dd_0^2}dt
\right .
$$
$$
\left . \qquad -2\int_0^{\infty}\sqrt{t}{\cdot}{\rm Tr} \ {\dot
\Dd}_0\Dd_0^2e^{-t\Dd_0^2}dt \right \}
$$
$$
={2\over{\sqrt{\pi}}}{\cdot}\int_0^{\infty}{d\over{dt}}\{\sqrt{t}{\cdot}{\rm Tr}
\ {\dot \Dd}_0e^{-t\Dd_0^2}\}dt
$$
$$
=-{2\over{\sqrt{\pi}}}{\cdot}\lim_{\e \to 0}\sqrt{\e}{\cdot}{\rm Tr} \ {\dot
\Dd}_0e^{-{\e}\Dd_0^2} \ \ .$$

\bigskip

Another formula for ${d\over{dr}}\{\eta_{\Dd_r}(0)\}|_{r=0}$ is the result
of the asymptotic expansion of
${\rm Tr} \ {\dot \Dd}_0e^{-t\Dd_0^2}$
(see (\ref{e:expansion1})).  Assume, for
instance, that ${\dot \Dd}$ is of order $1$, then

$$
 {\rm Tr} \ {\dot \Dd}_0e^{-t\Dd_0^2} = \sum_{k=0}^N
 t^{{k -n-1}\over{2}}b_k + O(t^{N - {\frac {n-1}{2}}}) \ \ .$$

\noi
We differentiate
$$
{d\over{dr}} \left
\{{1\over{\G({{{s+1}\over{2}})}}}\int_0^{\infty}t^{{s-1}\over{
2}} {\rm Tr} \ \Dd_re^{-t\Dd_r^2}dt \right \}|_{r=0}
$$
$$
= {1\over{\G({{{s+1}\over{2}})}}}\int_0^{\infty}t^{{s-1}\over{2}}
{\rm Tr} \ {\dot \Dd}_0e^{-t\Dd_0^2}dt
$$
$$
\qquad +
{2\over{\G({{{s+1}\over{2}})}}}\int_0^{\infty}t^{{s+1}\over{2}}
\frac{d}{dt}{\rm Tr} \ {\dot \Dd}_0e^{-t\Dd_0^2}\}dt
$$
$$
= {1\over{\G({{{s+1}\over{2}})}}}\int_0^{\infty}t^{{s-1}\over{2}}
{\rm Tr} \ {\dot \Dd}_0e^{-t\Dd_0^2}dt
$$
$$
+{2\over{\G({{s+1}\over{2}}})}(t^{\frac{s+1}{2}}{\rm Tr} \ {\dot
\Dd}_0e^{-t\Dd_0^2}]^{\infty}_0)
$$
$$
-{s+1\over{\G({{{s+1}\over{2}})}}}\int_0^{\infty}t^{{s-1}\over{2}}
{\rm Tr} \ {\dot \Dd}_0e^{-t\Dd_0^2}dt
$$
$$
={2\over{\G({{s+1}\over{2}}})}(t^{\frac{s-1}{2}}{\rm Tr} \ {\dot
\Dd}_0e^{-t\Dd_0^2}]^{\infty}_0)
$$
$$
-{{s}\over{\G({{{s+1}\over{2}})}}}\int_0^{\infty}t^{{s-1}\over{2}}{\rm
Tr} \ {\dot \Dd}_0e^{-t\Dd_0^2}dt \ \ .
$$

\bigskip

\noi The first term on the last line is equal to $0$ for $Re(s)$
large enough and does not affect the meromorphic extension of
${d\over{dr}}\{\eta_{\Dd_r}(s)\}|_{r=0}$. The second term gives us
what we need

$$
- \lim_{s \to 0} \ {{s}\over{\G({{{s+1}\over{2}})}}}\int_0^{\infty}{\rm Tr} \
t^{{s-1}\over{2}}{\dot \Dd}_0e^{-t\Dd_0^2}dt
$$
$$
=- \lim_{s \to 0} \
\frac{s}{\Gamma(\frac{s+1}{2})}\int_0^1t^{{s-1}\over{2}}
\sum_{k=0}^N t^{{k - n-1}\over{2}}b_k dt
$$
$$
=- {2\over{\sqrt{\pi}}}{\cdot}\lim_{s \to 0} \
s{\cdot}\sum_{k=0}^N{{b_k}\over{s  +k - n}} = -
{{2b_{{n}}\over{\sqrt{\pi}}}} \ .$$

\bigskip

\noi In particular the variation disappears if \noi $n = dim \ M$
is even by the theorem 1.13.2 in \cite{Gi95}. We get the same
result in the case of ${\dot \Dd}$ of order $0$, i.e.

$${d\over{dr}}\{\eta_{\Dd_r}(0)\}|_{r=0} = -
{{{2c_n}}\over{\sqrt{\pi}}} \ .$$

\noi
where $\{c_k\}$ is the set of new coefficients.

\bigskip

Now, let us discuss the last ingredient in the $\z$-determinant of the
Dirac operator - the modulus of
${\rm det}_{\z}\Dd$ - the (square root of the)
determinant of $\Dd^2$ . We have
already written the formula

$$
{\rm det}_{\z}\Dd^2 = e^{-\z_{\Dd^2}'(0)} \ \ .$$

\noi
Let us remind the reader that ${\frac d{ds}}\z_{\Dd^2}(s)|_{s=0}$ is given
by the formula

\begin{equation}\label{e:dt1}
{\frac d{ds}}\z_{\Dd^2}(s)|_{s=0} = \int_0^{\infty} {\frac
1{t}}{\rm Tr} \ e^{-t\Dd^2} dt \ \,
\end{equation}

\noi under the assumption $ker(\Dd)=0$ and $dim \ M$ is odd. Let
us explain how to interpret formula (\ref{e:dt1}). The trace ${\rm
Tr} \ e^{-t\Dd^2}$ has an asymptotic expansion given by
(\ref{e:expansion0}), which leads to a meromorphic extension of
the $\z$-function to the whole complex plane. Lemma \ref{l:z}
tells us that $\z_{\Dd^2}(s)$ is holomorphic in the neighborhood
of $s=0$ , hence the derivative with respect to $s$ exists. Let
$\k_{\Dd^2}(s)$ denote the integral $\int_0^{\infty}t^{s-1}{\rm
Tr} \ e^{-t\Dd^2}dt$. The formula (\ref{e:zreg}) gives us the
expansion of $\k_{\Dd^2}(s)$ in the neighborhood of $s = 0$. We
have

$$\k_{\Dd^2}(s)= {{a_{{n\over{2}}}}\over{s}} + h_1(s) \ \ .$$

\noi
Now, the derivative of the $\z$-function at $s=0$ is obtained as follows:

$$\z_{\Dd^2}'(0) = {\frac d{ds}} {\frac{\k_{\Dd^2}(s)}{\G(s)}}|_{s=0}
$$

$$
={\frac d{ds}}({\frac{a_{n\over{2}} + s(\k_{\Dd^2}(s) -{\frac
{a_{n\over{2}}}{s}})}{1 +
s\g + s^2h(s)}})|_{s=0}
$$
$$=(\k_{\Dd^2}(s) -{\frac{a_{n\over{2}}}{s}})|_{s=0} - {\g}a_{n\over{2}} .
$$

\bigskip

\noi
If $n$ is odd then the coefficient $a_{n\over{2}} = 0$ and we can
(``formally") write

\begin{equation}\label{e:formal}
- {\rm ln} \; {\rm det}_{\z}\Dd^2 = \k_{\Dd^2}(s)|_{s=0} = \int_0^{\infty}
{\frac 1{t}}{\rm Tr} \ e^{-t\Dd^2} dt.
\end{equation}

\bigskip
It is worth mentioning that the variation of ${\rm det}_{\z}\Dd^2$
is by no means a local invariant. Assume that we have a family of
invertible Dirac operators $\{\Dd_r\}$ , then we can use Duhamel's
Principle as in the case of the $\eta$-invariant. We obtain

$$
\frac d{dr}\{{\rm ln} \ {\rm det}_{\z}\Dd_r^2\}|_{r=0} = - \frac
d{dr} \left . \int_0^{\infty} {\frac 1{t}}{\rm Tr} \ e^{-t\Dd_r^2}
dt \right |_{r=0}
$$
$$
=2\int_0^{\infty}{\rm Tr} \ {\dot \Dd}_0\Dd_0 e^{-t\Dd_0^2}\
dt\qquad\qquad
$$
$$
= -2\int_0^{\infty}\frac d{dt} \ {\rm Tr} \ {\dot
\Dd}_0\Dd_0^{-1}\{e^{-t\Dd_0^2}\}\ dt \ .$$

\bigskip

\noi This gives us the formula

\begin{equation}\label{e:et1}
\frac d{dr}\{{\rm ln} \ {\rm det}_{\z}\Dd^2\}|_{r=0} = 2{\cdot}\lim_{\e \to
0} {\rm Tr} \ {\dot \Dd}_0\Dd_0^{-1}e^{-\e \Dd_0^2}.
\end{equation}

\noi This formula allows us to see that ${\rm det}_{\z}\Dd^2$ is
actually a highly non-local invariant as it involves the kernel of
the operator $\Dd_0^{-1}$ .

\bigskip
To give a simple example let us consider the family
$\{\Delta_r = \Dd^2 e^{r\a}\}_{0 \le r \le 1}$, where
$\a: C^{\infty}(M;S) \to C^{\infty}(M;S)$ is an
operator with smooth kernel. We repeat the computations which lead
to (\ref{e:et1}) and obtain

$$\frac d{dr}\{{\rm ln} \; {\rm det}_{\z}\Delta_r\}
= {\rm Tr} \a \ \ .$$

\bigskip

\noi which implies
$$
{\rm ln} \ {\rm det}_{\z}\Delta_1 - {\rm ln} \ {\rm det}_{\z}\Delta_0 =
\int_0^1 {\rm Tr} \ \a dr = {\rm Tr} \ \a .
$$

\noi We have proved the equality

\begin{equation}\label{e:et2}
{\rm det}_{\z}\Dd^2 e^{\a}
= {\rm det}_{\z}\Dd^2{\cdot}{\rm det}_{Fr}e^{\a}.
\end{equation}

\bigskip
On the other hand let us discuss the $\eta$-invariant for the family
$\{\Dd_r = \Dd + r\a\}$ . We have

$${d\over{dr}}\{\eta_{\Dd_r}(0)\} = -{2\over{\sqrt{\pi}}}{\cdot}\lim_{\e
\to 0}\sqrt{\e}{\cdot}{\rm Tr} \ {\dot \Dd}_0e^{-{\e}\Dd_0^2}$$

$$= -{2\over{\sqrt{\pi}}}{\cdot}\lim_{\e \to 0}\sqrt{\e}{\cdot}{\rm Tr} \ \a
e^{-{\e}\Dd_0^2}$$

$$ = -{2\over{\sqrt{\pi}}}{\cdot}\lim_{\e \to 0}\sqrt{\e}{\cdot}
{\rm Tr} \ \a = 0
\ ,$$

\noi
and as a result

\begin{equation}\label{e:et3}
\ \ \ \ \ \ \ \ \ \ \ \ \
\eta_{\Dd + \a}(0) = \eta_{\Dd}(0).
\end{equation}

\bigskip

\section{Determinants of Dirac operators on a manifold with
boundary}\label{s:bd}

\bigskip

In this Section we discuss the determinants of Dirac operators on
a manifold with boundary. The new ingredient is that, in order to
get a nice elliptic operator out of $\Dd$ , we have to consider
the boundary conditions. The choice of boundary condition
determines the domain of the operator $\Dd$ . We will not discuss
here the most general space of elliptic, self-adjoint boundary
conditions for $\Dd$ introduced in the recent work of Kirk and
Lesch (see \cite{KL00}). We stick to the more conventional
Grassmannian of the boundary conditions of Atiyah--Patodi--Singer
type. We avoid also a discussion of the case of non-product metric
structures in the neighborhood of the boundary, which rises to the
table many unpleasant analytical issues.

An unexpected advantage of the fact that we discuss boundary problems is that
in our situation $det_{\z}$ is
in fact equal (up to a scalar)  to the true Fredholm determinant.

Let $M$ denote an odd-dimensional compact manifold with boundary
$Y$ and $\Dd: C^{\infty}(M;S) \to C^{\infty}(M;S)$ a compatible
Dirac operator acting on sections of $S$ , a bundle of Clifford
modules over $M$ . Assume that the Riemannian metric on $M$ and
the Hermitian structure on $S$ are products in a certain collar
neighborhood of the boundary. Let us fix a parameterization $N
=[0,1] \times Y$ of the collar. Then, in $N$, the operator $\Dd$
has the form

\begin{equation}\label{e:Dd}
\ \ \ \ \ \ \ \ \ \ \ \ \ \ \
\Dd = G(\partial_u + B) \ \, ,
\end{equation}

\noi
where $G : S|Y \to S|Y$ is a unitary bundle isomorphism (Clifford
multiplication by the unit normal vector) and
$B: C^{\infty}(Y;S|Y) \to  C^{\infty}(Y;S|Y)$ is the corresponding Dirac
operator on $Y$, an elliptic self-adjoint operator of first order.
Furthermore, $G$ and $B$
do not depend on the normal coordinate $u$ and they satisfy the identities

\begin{equation}\label{e:GB}
\ \ \ \ \ \ \
G^2 = -Id \ \ \ {\rm and} \ \ \ GB = -BG \ \, .
\end{equation}

\noi
Since $Y$ has dimension $2m$ the bundle $S|Y$ decomposes into its positive
and negative chirality components
$S|Y = S^+ \bigoplus S^-$ and we have a corresponding splitting of the
operator $B$ into
$B^{\pm}: C^{\infty}(Y;S^{\pm}) \to C^{\infty}(Y;S^{\mp})$ , where
$(B^+)^* = B^-$. The operator (\ref{e:Dd})
can be rewritten in the form

\begin{equation}\label{e:ch1}
\pmatrix{
i & 0 \cr
0 & -i \cr}
\left(
\partial_u +
\pmatrix{ 0 & B^- \cr B^+ & 0 \cr} \right). \ \,
\end{equation}

In order to obtain an unbounded Fredholm operator with sufficient
regularity properties we have to impose
a boundary condition on the operator $\Dd$ . Let $\Pi_>$ denote the spectral
projection of $B$ onto the subspace of $L^2(Y;S|Y)$
spanned by the eigenvectors corresponding to the nonnegative eigenvalues of
$B$. It is well known that $\Pi_>$ is
an elliptic boundary condition for the operator $\Dd$ (see \cite{AtPaSi75},
\cite{BoWo93}). The meaning of ellipticity is described below. We introduce
the unbounded operator $\Dd_{\Pi_>}$ equal to the operator $\Dd$ with domain

\[ {\rm dom} \ \Dd_{\Pi_>}
= \{s \in H^1(M;S) \ ; \ \Pi_>(s|Y) = 0 \} \,, \]

\noi
where $H^1$ denotes the first Sobolev space. Then the operator

$$
\Dd_{\Pi_>} = \Dd : {\rm dom}(\Dd_{\Pi_>}) \to L^2(M;S)$$

\noi
is a Fredholm operator with kernel and cokernel consisting only of smooth
sections.

The orthogonal projection $\Pi_>$ is a pseudodifferential operator
of order 0 (see \cite{BoWo93}). Let us point out that  we can take
any pseudodifferential operator $R$ of order $0$ with principal
symbol equal to the principal symbol of $\Pi_{>}$ and obtain an
operator $\Dd_R$ which satisfies the aforementioned properties. In
the following, however, we concentrate on the specific subset of
the space of self-adjoint elliptic boundary conditions. There
exists another pseudodifferential projection on $Y$ , which is in
fact the central object in the theory of  elliptic boundary value
problems. Let us briefly explain this point. In contrast to the
case of an elliptic operator on a closed manifold, the operator
$\Dd$ has an infinite-dimensional space of solutions. More
precisely, the space

$$\{s \in C^{\infty}(M:S) \ ; \ {\Dd}s = 0 \ \ in \ M \setminus Y\}$$

\noi is infinite-dimensional. We introduce the Calde-ron
projection, which is the projection onto $\Hh(\Dd)$ of the {\it
Cauchy Data space} of the operator $\Dd$

\bigskip

$$\Hh(\Dd) = \{f \in C^{\infty}(Y;S|Y) \ ; \ \exists \ {s \in
C^{\infty}(M;S)} $$
$$\qquad \rm{s.t.} \ \ \Dd(s) = 0 \ {\rm in} \
M \setminus Y \ {\rm and} \ s|Y = f\} \ .$$

\bigskip

\noi The projection $P(\Dd)$ is a pseudodifferential operator with
principal symbol equal to the symbol of $\Pi_>$ . It is also an
orthogonal projection in the case of a Dirac operator on an
odd-dimensional manifold (see \cite{BoWo93}). The operator $\Dd$
has the {\it Unique Continuation Property}, and hence we have an
one to one correspondence between solutions of the operator $\Dd$
and the traces of solutions on the boundary $Y$. This roughly
explains why only the projection $\Pp_R$ onto the kernel of the
boundary conditions $R$ matters. If the difference $\Pp_R -
P(\Dd)$ \hskip 2mm is an operator of order $-1$ , then it follows,
that by choosing the domain of the operator $\Dd_R$ as above, we
throw away almost all solutions of the operator $\Dd$ on $M
\setminus Y$, with the possible exception of a finite dimensional
subspace. The above condition on $\Pp_R$ also allows us to
construct a parametrix for the operator $\Dd_R$ , hence we obtain
regularity of the solutions of the operator $\Dd_R$. We refer to
\cite{BoWo93} for more details.

This explains why in \cite{SSKPW299} we restricted ourselves to
the study of the Grassmannian $Gr^*_{\infty}(\Dd)$ of all
orthogonal pseudodifferential projections $P$ such that

\begin{equation}\label{e:grsa}
\qquad P - P(\Dd) \ \ {\rm is}  \ {\rm a} \ {\rm smoothing} \ {\rm
operator} \ \,
\end{equation}

$\qquad \qquad {\rm and}  \ \ \ - GPG = Id - P \ \ .$

\bigskip
The first condition implies the ellipticity of the operator $\Dd_P$
and the second guarantees self-adjointness.
The spectral projection $\Pi_>$ is an element of $Gr^*_{\infty}(\Dd)$ if
and only if $ker \ B = \{0\}$.

\bigskip

\begin{rem}\label{r:ch1}
Again let us point out that the space $Gr_{\infty}^*(\Dd)$ is far
from being the space of all elliptic boundary conditions for the
Dirac operator $\Dd$ . An important example is given by the
condition determined by chirality (see (\ref{e:ch1})). The
operator $P_{\pm} = {\frac 1{2}}(Id \mp i\G)$ is the orthogonal
projection of $S|_Y$ onto $S^{\pm}$ and provides $\Dd$ with a
(local) {\it chiral} elliptic boundary condition. This means that
the operator $\Dd_{\pm} = \Dd$ with domain
$$
{\rm dom} \ \Dd_{\pm} = \{s \in H^1(M;S) \mid P_{\pm}(s|_Y) = 0\} \,,$$

\noi
is Fredholm and that its kernel and cokernel consist of smooth sections
only. The operators $\Dd_{\pm}$ are not self-adjoint,
but we have the equalities

\begin{equation}\label{e:dpm}
\Dd_+^* = \Dd_-  \ \ {\rm and} \ \ {\rm index} \ \Dd_{\pm} = 0. \ \,
\end{equation}

\noi
It is not difficult to see that $\Delta_{\pm} = \Dd_{\mp}\Dd_{\pm}$ is
equal to the operator
$\Dd^2$ with Dirichlet (resp. Neumann) condition on $S^+$ and Neumann
(resp. Dirichlet) condition on $S^-$.
\end{rem}

\bigskip

For any $P \in Gr^*_{\infty}(\Dd)$ the operator $\Dd_P$ has a discrete
spectrum nicely distributed along
the real line. It was shown by the second author that $\eta_{\Dd_P}(s)$
and $\z_{\Dd_P^2}(s)$ are well-defined functions, holomorphic for
$Re(s)$ large and having meromorphic extensions to the whole complex
plane with only simple poles. In particular
both functions are holomorphic in a neighborhood of $s = 0$. Therefore
${\rm det}_{\z}\Dd_P$ is a well-defined, smooth function on $Gr_{\infty}^*(\Dd)$
(see \cite{KPW99}). We will discuss the regularity of $\eta$-function of
the operator
$\Dd_P$ , with $P \in Gr_{\infty}^*(\Dd)$ in Section 9.
Now we discuss  the ``true"
determinant, which lives on the space $Gr_{\infty}^*(\Dd)$ .

\bigskip

The determinant line bundle over the space of Fredholm operators
was first introduced in a seminal paper of Quillen \cite{Qu85}. An
equivalent better suited to our purposes was subsequently given
by Segal (see \cite{Seg90}), and we follow his approach.
Let ${\rm Fred}(\Hh)$ denote the space of Fredholm
operators on a separable Hilbert space $\Hh$.  First we work in
the connected component ${\rm Fred}_{0}(\Hh)$ of this space
parameterizing operators of $index$ zero. For $A\in {\rm Fred}_{0}(\Hh)$
define

$$
{\rm Fred}_A = \{S \in {\rm Fred}(\Hh) \ ; \
S - A \ \ {\rm is} \ {\rm trace-class} \} \
\ .$$

\bigskip

\noi
Fix a trace-class operator $\Aa$ such that $S = A + \Aa$ is
an invertible operator. Then the determinant line of $A$ is defined as

\begin{equation}\label{e:1.1}
{\rm Det} \ A = {\rm Fred}_A \times \C/_{\cong} \ \,
\end{equation}

\bigskip

\noi
where the equivalence relation is defined by

$$(R , z) = ((RS^{-1})S , z) \simeq (S , z{\cdot}
{\rm det}_{Fr}(RS^{-1})) \ \ .$$

\bigskip

\noi
The Fredholm determinant of the operator $RS^{-1}$ is
well-defined, as it is of the form $Id_{\Hh}$ plus a trace class
operator. Denoting the equivalence class of a pair $(R,z)$ by
$[R,z]$, complex multiplication is defined on ${\rm Det} \ A$ by

\begin{equation}\label{e:multn}
   \lambda{\cdot}[R,z] = [R,\lambda z].
\end{equation}

\noi
The {\em canonical determinant element} is defined by

\begin{equation}\label{e:detelt}
{\rm det} \ A := [A,1] ,
\end{equation}

\noi
and is non-zero if and only if $A$ is invertible. The complex lines fit
together over
${\rm Fred}_{0}(\Hh)$ to define a
complex line bundle $\Ll$, the determinant line bundle. To see
this, observe first that over the open set $U_{\Aa}$ in
${\rm Fred}_{0}(\Hh)$ defined by

$$
U_{\Aa} = \{F \in {\rm Fred}_{0}(\Hh) \ ; \ F + \Aa \ \ {\rm is} \
{\rm invertible} \} ,$$

\bigskip

\noi
the assignment \hskip 2mm $F \to {\rm det} \ F$ \hskip 2mm defines
a trivializing (non-vanishing) section of $\Ll_{|U_{\Aa}}$. The
transition map between the canonical determinant elements over
$U_{\Aa} \cap U_{\Bb}$ is the smooth (holomorphic) function

$$g_{{\Aa}{\Bb}}(F) = {\rm det}_{Fr}((F+\Aa)(F+\Bb)^{-1}) \ \ .$$

\bigskip

\noi
This defines $\Ll$ globally as a complex line bundle over
${\rm Fred}_{0}(\Hh),$ endowed with the canonical section \hskip 2mm $A
\to {\rm det} \ A$ . If  ${\rm ind} \ A = d$
we define ${\rm Det} \ A $ to be the
determinant line of  $A \oplus 0$ as an operator $\Hh\too\Hh\oplus
\CC^{d}$ if  $d > 0$ , or  $\Hh \oplus \CC^{-d}\too\Hh$ if  $d <
0$ and the construction extends in the obvious way to the other
components of ${\rm Fred} (\Hh)$. Note that the canonical section is
zero outside of ${\rm Fred}_{0}(\Hh).$

\bigskip

We use this construction in order to define the determinant line
bundle over $Gr_{\infty}(\Dd)$. For each projection $P \in
Gr_{\infty}(\Dd)$ we have the (Segal) determinant line ${\rm Det}(P(\Dd),P)$ of
the operator

$$\Ss(P) = PP(\Dd) : \Hh(\Dd) \to {\rm Ran} \ P \ \ ,$$

\noi
and the determinant line ${\rm Det}\ \Dd_P$ of the boundary-value problem
$\Dd_P: {\rm dom} \ (\Dd_P) \too L^{2}(M;S)$. These lines fit together
in the manner explained above to define determinant line bundles
${\rm DET}_{P(\Dd)}$ and ${\rm DET} \ \Dd$, respectively, over the
Grassmannian (some care has to be taken as the operator acts
between two different Hilbert spaces, but with the obvious
notational modifications we once again obtain well-defined
determinant line bundles). The topology of the Grassmannians
(see \cite{BoWo89}, \cite{DoWo91}) implies that the bundle
${\rm DET}_{P(\Dd)}$ is a non-trivial line bundle over
$Gr_{\infty}(\Dd)$, but when restricted to the Grassmannian
$Gr_{\infty}^*(\Dd)$ it is canonically trivial.  The canonical section becomes
a function in this trivialization. We call this function {\it Canonical
Determinant} and we denote
its value at $P$ by ${\rm det}_{\Cc}\Dd_P$.

Now we give more precise description of ${\rm det}_{\Cc}\Dd_P$. Simon Scott
showed that elements of $Gr^*_{\infty}(\Dd)$ are
in one to one correspondence with the unitary elliptic operators $T :
C^{\infty}(Y;S^+) \to C^{\infty}(Y;S^-)$ ,
which satisfy an additional condition (see \cite{Sc95}).  Namely,
let us introduce the operator $V_> = (B^+B^-)^{-1}B^+$ . We assume that

$$T - V_> \ \ \ {\rm is} \ {\rm a} \ {\rm smoothing} \ {\rm operator}.
$$

\noi
The correspondence is as follows: if we fix the operator $T$ as above
then the corresponding projection is

\begin{equation}\label{e:corr}
T \to P = {1\over{2}}
\pmatrix{
Id_{F^+} & T^{-1} \cr
T & Id_{F^-} \cr}
\ \, .
\end{equation}

\noi
Let us stress that the invertibility assumption on the tangential
operator $B$ can be easily relaxed when we discuss this construction
(see Section 7.3 of \cite{SSKPW299} for the details). Let us also point
out that this fixes the isomorphism of
$Gr^*_{\infty}(\Dd)$ with $U^{\infty}(F^-)$ the group of unitary operators
on the sections of $S^-$ of the form \hskip 2mm $Id_{F^-}$ plus {\it
smoothing operator} .
Let $K : C^{\infty}(Y;S^+) \to C^{\infty}(Y;S^-)$ be a unitary
operator such that $\Hh(\Dd) = graph \ K$ . Then the operator

\begin{equation}\label{e:tr}
U(P) = \pmatrix{
Id_{F^+} & 0 \cr
0 & TK^{-1} \cr}
\ \,
\end{equation}

\bigskip
\noi
has the property

$$P = U(P)P(\Dd)U(P)^{-1}$$

\noi
and it defines an isomorphism \hskip 3mm $P \to TK^{-1}$ \hskip 3mm between
$Gr^*_{\infty}(\Dd)$ and $U^{\infty}(F^-)$ . Now we have a
well-defined operator

$$U(P)^{-1}\Ss(P): \Hh(\Dd) \to \Hh(\Dd) \ \ .$$

\noi It is of the form \hskip 2mm $Id_{\Hh(\Dd)}$ plus {\it
smoothing operator}, hence it has a well-defined Fredholm
determinant and straightforward computations show that

$$
{\rm det}_{Fr}U(P)^{-1}\Ss(P) = {\rm det}_{Fr}
\left({{Id + KT^{-1}}\over{2}}\right) .
$$

\noi
All this was explained in Section 1 of \cite{SSKPW299}. The study of the
preferred trivialization, defined by means of
the operator $U(P)$ , now shows that we have the equality

\begin{equation}\label{e:det1}
{\rm det}_{\Cc}\Dd_P = {\rm det}_{Fr}U(P)^{-1}\Ss(P) \ \, .
\end{equation}

\bigskip
The question arises: Is  $det_{\Cc}$ related to $det_{\z}$? A
positive answer was given in work of Scott and Wojciechowski, as
the main result of \cite{SSKPW299} is

\bigskip

\begin{thm}\label{t:0.1}
The following equality holds over $Gr_{\infty}^*(\Dd)$:

\begin{equation}\label{e:t1}
{\rm det}_{\z}\Dd_P = {\rm det}_{\z}\Dd_{P(\Dd)}
{\cdot}{\rm det}_{\Cc}\Dd_P \ \, .
\end{equation}

\end{thm}

\bigskip

To prove Theorem \ref{t:0.1} we study the variation of the determinants.
More precisely,
we fix two projections $P_1 , P_2 \in Gr_{\infty}^*(\Dd)$ such that the
operators $\Dd_{P_i}$
are invertible. Next, we choose a family of unitary operators of the form

$$\left \{ \pmatrix{
Id_{F^+} & 0 \cr
0 & g_r \cr} \right
\}_{0 \le r \le 1}
\ \ ,$$

\noi
where $g_r : F^- \to F^-$ is a unitary operator, and such
that $g_r - Id_{F^-}$ is an operator with a smooth kernel for any
$r$ , and $g_0 = Id_{F^-}$ . We define two families of boundary
conditions:

$$P_{i,r} =
\pmatrix{
Id_{F^+} & 0 \cr
0 & g_r \cr}
P_i
\pmatrix{
Id_{F^+} & 0 \cr
0 & g_r^{-1} \cr}
\ \ ,$$

\noi
and study the relative variation:

\begin{equation}\label{e:rv}
{d\over{dr}}\{{\rm ln} \ {\rm det} \ \Dd_{P_{1,r}}
- {\rm ln} \ {\rm det} \
\Dd_{P_{2,r}}\}|_{r=0}
\ \,
\end{equation}

\noi
for both the {\it Canonical determinant} and the
$\z$-determinant. Of course we face the technical problem of
dealing with a family of unbounded operators with varying domain.
To circumvent this, and to make sense of the variation with respect to
the boundary condition we follow Douglas and Wojciechowski
\cite{DoWo91} and apply their {\it ``Unitary Trick"}. It is not difficult to
define an extension of our family of unitary operators on the boundary
sections to a family $\{U_r\}$  of unitary operators acting on
$L^2(M;S)$ (see Section 9 for more details). The operator
$\Dd_{P_{i,r}}$ is unitarily equivalent to the operator
$(\Dd_r)_{P_i}$ , where

$$\Dd_r = U_r^{-1}{\Dd}U_r \ \ .$$

\bigskip
\noi
Both the $\z$-determinant and the canonical determinant are
invariant under this unitary twist which allows us to show that
both determinants have variation given by the same expression

\begin{equation}\label{e:vr}
{d\over{dr}}\{{\rm ln} \ {\rm det} \ \Dd_{P_{1,r}}
- {\rm ln} \ {\rm det} \
\Dd_{P_{2,r}}\}|_{r=0}
\ \,
\end{equation}

$= {\rm Tr} \ {\dot D}_0 (\Dd_{P_1}^{-1} - \Dd_{P_2}^{-1}) \ \ ,$

\bigskip
\noi where ${\dot \Dd}_0$ denotes the operator
${d\over{dr}}\Dd_r|_{r=0}$. Now we use the fact that the set of
projections $P \in Gr_{\infty}^*(\Dd)$, for which the operator
$\Dd_P$ is invertible is actually path connected (see Section 7.2
of \cite{SSKPW299}) and integrate the equality

$${d\over{dr}}\{{\rm ln} \ {\rm det}_{\z} \ \Dd_{P_{1,r}}
- {\rm ln} \ {\rm det}_{\z} \
\Dd_{P_{2,r}}\}|_{r=0} $$

$= {d\over{dr}}\{{\rm ln} \ {\rm det}_{\Cc} \
\Dd_{P_{1,r}}
- {\rm ln} \ {\rm det}_{\Cc} \ \Dd_{P_{2,r}}\}|_{r=0} \ \ ,$

\bigskip

\noi
in order to obtain formula (\ref{e:t1}) of Theorem \ref{t:0.1}.

The reader might think that formula (\ref{e:vr}) is incorrect as the
variation of the phase of
the $\z$-determinant is not present. However, we will see in Section 9 that
the variation of the $\eta$-invariant in our situation does depend only on
$\{g_r\}$, and
not on the choice of the base projection, hence the variation here is the
same at $P_1$ as it is at $P_2$.
We learn more about the properties of the $\eta$-invariant on the
Grassmannian in Section \ref{s:eta}.

\bigskip

\section{An outline of the method}\label{s:m}

\bigskip

The idea to use the adiabatic limit in this particular
way belongs to Singer (see \cite{Si88}).
There are three basic ingredients which we use in our approach to
the decomposition of the $\z$-determinant.

\bigskip

First, we rely heavily on the assumption that metric structures are product near
the boundary. This implies that the operator $\Dd$ has a cylindrical form in
the collar neighborhood
of the boundary. The determinant is expressed via different {\it Heat Operators}
determined by $\Dd$ and those operators are not local. The crucial quantity
here is

$$\int_M {\rm tr} \ \Ee(t;x,x)dx \ \ ,$$

\noi
where $\Ee(t;x,y)$ denotes the kernel of such an operator. We know the
construction of $\Ee(t;x,y)$ on
a closed manifold, hence in the interior of $M$. The product structure
gives also the explicit formulas for
the kernel $\Ee(t;x,y)$ on the cylinder. The problem is to paste those
kernels in order to get a kernel on $M$. Moreover, the
endomorphism $\Ee(t;x,y)$ is not determined
via coefficients of $\Dd$ at $x$ and $y$ only but depends on global
information from the whole manifold $M$. Now, the construction
of the kernel $\Ee(t;x,y)$ on a closed
manifold $M$ is standard and described in many different places.
What is important for us is that the estimates, obvious in the case
of flat space, hold also in the case of a general manifold.

\bigskip

\begin{prop}\label{p:est}
Let $\Dd$ be a Dirac operator on a closed manifold and $\Ee(t;x,y)$ and
$\Ff(t;x,y)$ denote the kernels
of the operators $e^{-t\Dd^2}$ and $\Dd e^{-t\Dd^2}$ . Then there exist
positive constants $c_1$ and $c_2$
such that

$$\|\Ee(t;x,y)\| \le c_1t^{-{n\over{2}}}e^{-c_2{{d^2(x,y)}\over{t}}},$$

$$\|\Ff(t;x,y)\| \le c_1t^{-{{n+1}\over{2}}}e^{-c_2{{d^2(x,y)}\over{t}}}$$

\bigskip

\noi
for any $x,y \in M$ and any $t > 0$ .

\end{prop}

\bigskip

We show that those estimates extend easily to our situation. We refer to
\cite{BoWo93}, \cite{DoWo91} and \cite{Ta296} for additional information
on this subject and more comprehensive bibliography as the literature
on this topic is extremely rich.

\bigskip

Second, we use Duhamel's Principle to paste kernels. The Duhamel
Principle shows explicitly that the heat kernel on $M$ splits into
interior part, cylindrical part and the error term.
It also provides us the tools to study the error term.

\bigskip

Third, we assume that the tangential operator $B$ is invertible.
This allows us to make assumptions concerning the behavior of the
eigenvalues of the boundary problems, which eventually allows us
to discard the {\it large time contribution}. We also rule out the
existence of the $L^2$-solutions of $\Dd$ on manifolds with
cylindrical ends, which enter the picture during our analysis.
These assumptions secure the non-existence of the ``small"
eigenvalues in the situations we study.

\bigskip
We will not discuss here the scattering theory, which enters the
picture in the case of the non-invertible tangential operator. The
situation is as follows. Assume that we have a given decomposition
of a closed manifold $M$ into two submanifolds $M_1$ and $M_2$
along the submanifold $Y$ of codimension $1$. As we stretch the
collar neighborhood around $Y$, the eigenvalues of the Dirac
operator $\Dd$ change. The assumptions we made above guarantee
that they stay bounded away from $0$ . However if the tangential
operator is non-invertible, we have to deal with small eigenvalues
of $\Dd$. They fall in two different categories. We have finitely
many eigenvalues which decay exponentially with respect to the
length $R$ of the cylinder $[-R,R] \times Y$ joining $M_1$ and
$M_2$. These eigenvalues are constructed from $L^2$ solutions of
$\Dd$ on $M_1$ and $M_2$ with the cylinders of infinite length
attached. There are also an infinite family of eigenvalues of size
$\frac 1{R}$, which can be constructed from the eigenvalues of the
corresponding Dirac operators on the circles of large radius
determined by the {\it Scattering Theory} defined by $\Dd$. This
type of analysis was applied to the Atiyah--Patodi--Singer problem
by Werner M$\ddot {\rm u}$ller (see \cite{Mu94}, see also
\cite{Mu98} for related results). Following M$\ddot {\rm u}$ller
the authors were able to establish a decomposition formula for the
$\z$-determinant in the case of non-invertible tangential
operator. We refer to the recent paper \cite{JKPW3} (see also
\cite{JKPW4}) for more details.

\bigskip

\section{Cylinder, Duhamel's Principle and Heat Kernels on a manifold with
boundary}\label{s:cyl}

We start with the infinite cylinder $[0,\infty) \times Y$ and the operator
$\Dd^2 = -\partial_u^2 + B^2$ ,
subject to the boundary condition at $u = 0$ . We collect several
explicit formulas for the kernels of
the heat operators  determined by $\Dd^2$ subject to different boundary
conditions. Then we show how Duhamel's  Principle leads to
the splitting of the trace of the heat operator on a manifold $M$ onto
interior contribution,  cylinder contribution and the error term.

We start with the Dirichlet condition. We introduce the operator $\Delta_d
= \Dd^2$ with the domain

$$
{\rm dom} \  \Delta_d = \{s \in C^{\infty}([0,\infty) \times Y;S) ; s|_{u=0} =
0\} \ \ .$$

\noi
It has a unique closed self-adjoint extension and therefore
$e^{-t\Delta_d}$ is well-defined and its kernel is given
by the formula

\begin{equation}\label{e:dd1}
\Ee_d(t;(u,x),(v,y)) = \end{equation}
$$ \qquad
{\frac1{\sqrt{4{\pi}t}}}\{e^{-\frac{(u-v)^2}{4t}}
 - e^{-\frac{(u+v)^2}{4t}}\}\ e^{-tB^2}(t;x,y) ,
$$

\bigskip

\noi where $e^{-tB^2}(t;x,y)$ denotes the kernel of the operator
$e^{-tB^2}$. Similarly, to discuss the Neumann condition we
introduce $\Delta_n = \Dd^2$ with domain

$$
{\rm dom} \  \Delta_n = \{s \in C^{\infty}([0,\infty) \times Y;S) ;
(\partial_us)|_{u=0} = 0\} .
$$

\noi The corresponding heat kernel is given by the formula
\begin{equation}\label{e:nn1}
\Ee_n(t;(u,x),(v,y)) =
\end{equation}
$$ \quad \frac1{\sqrt{4{\pi}t}}\{e^{-\frac{(u-v)^2}{4t}}
  + e^{-\frac{(u+v)^2}{4t}}\}\ e^{-tB^2}(t;x,y)
\ \ .$$

\bigskip

\noi Finally let us note the formula for the kernel of the
operators $e^{-t\Delta_{\pm}}$ (see Remark \ref{r:ch1})

$$
\Ee_{\pm}(t;(u,x),(v,y)) = \qquad\qquad\qquad\qquad\qquad
$$
$${\frac1{\sqrt{4{\pi}t}}}\{e^{-\frac{(u-v)^2}{4t}} \mp
e^{-\frac{(u+v)^2}{4t}}\}e^{-tB^2}(t;x,y)P_+ $$
$$+ {\frac1{\sqrt{4{\pi}t}}}\{e^{-\frac{(u-v)^2}{4t}} \pm
e^{-\frac{(u+v)^2}{4t}}\}e^{-tB^2}(t;x,y)P_- \ .$$

\bigskip

 The formulas for the Atiyah--Patodi--Singer condition are more
complicated. It follows from (\ref{e:ch1}) that the operator $B$
has a symmetric spectrum. Let $\{\m_n\}_{n \in \N}$ denote the set
of positive eigenvalues and $\{\phi_n\}$ the set of corresponding
eigenspinors, then the negative eigenvalues are $\{-\m_n\}$ with
the corresponding eigenspinors $\{G\phi_n\}$ . The heat kernel of
$\Dd_{\Pi_>}^2$ on the cylinder has the form

$$\sum_{n \in \N}g_n(t;u,v)\phi_n(x)\tensor \phi_n(y) $$
$$
+\sum_{n \in \N}g_{-n}(t;u,v)G\phi_n(x)\tensor G\phi_n(y) \ \ .$$

\bigskip

\noi
Recall the formulas for the functions $g_n(t;u,v)$ (see for instance
\cite{BoWo93}, (22.33) and (22.35))

\begin{equation}\label{e:g1}
g_n(t;u,v) =
{{e^{-{\m}_n^2t}}\over{2\sqrt{{\pi}t}}}{\cdot}\{e^{-{{(u-v)^2}\over{4t}}} -
e^{-{{(u+v)^2}\over{4t}}}\}
\ \,
\end{equation}

$ \ {\rm for} \ n>0$, and

\[
g_n(t;u,v) =
{{e^{-(-{\m}_n)^2t}\over{2\sqrt{{\pi}t}}}}{\cdot}\{e^{-{{(u-v)^2}\over{4t}}}
+
e^{-{{(u+v)^2}\over{4t}}}\} +
\ \,
\]

\[
(-\m_n)e^{-(-\m_n)(u+v)}{\cdot}{\rm erfc} \left({{u+v}\over{2\sqrt{t}}} -
(-\m_n)\sqrt{t} \right)
\ \,
\]

$ \ {\rm for} \ n<0 \ \ $ where

\[
{\rm erfc}(x) = {2\over{\sqrt{\pi}}}\int_x^{\infty} e^{-r^2}dr <
{2\over{\sqrt{\pi}}}e^{-x^2}
\ \, .
\]

\bigskip

\noi Note that all those kernels satisfy the estimates from
Proposition \ref{p:est}.

\bigskip

The interior heat kernel is defined by the kernel of the double of the
Dirac operator $\Dd$.
Let us recall that this operator has a natural double ${\tilde \Dd}$,
which leaves on $\tilde M$,
the double of a manifold $M$. Let ${\tilde \Ee}(t;x,y)$ denote the kernel
of the operator
$e^{-t{\tilde \Dd}^2}$ and let $\Ee_{cyl}(t;x,y)$ denote one of the kernels
on the cylinder
discussed above. We use them to construct the kernel of the operator
$e^{-t\Dd_{\Pi_.}^2}$
(and  the operators $e^{-t\Delta_d}$, $e^{-t\Delta_n}$,
$e^{-t\Delta_{\pm}}$) on
the manifold $M$. Roughly speaking we glue cylinder kernel  and interior
kernel together.

We introduce a smooth, increasing function
$\rho(a,b): [0, \infty) \to [0,1]$ equal to $0$ for $0 \le u \le a$
and equal to $1$ for $b \le u$ . We use $\rho(a,b)(u)$ to define

$$\phi_1 = 1 - \rho \left(\frac{5}{7} , \frac{6}{7}\right)  ,
\;  \psi_1 = 1
- \psi_{2},
$$

and

$$\phi_2 = \rho \left(\frac{1}{7} , \frac{2}{7}\right) , \; \psi_2 =
\rho \left(\frac{3}{7} , \frac{4}{7}\right).$$

\bigskip

\noi
We extend those functions to the symmetric functions on the whole
real line. All those functions are constant
outside the interval $[-1 , 1]$ and we use them to define the
corresponding functions on a manifold $M$.
Now we define $Q(t;x,y)$, a {\it ``Parametrix"} for the real heat
kernel ${\Ee}(t;x,y)$, by

\begin{equation}\label{e:p1}
Q(t;x,y) = \phi_1(x){\Ee}_{cyl}(t;x,y)\psi_1(y)
\end{equation}

$$+ \phi_2(x){\tilde \Ee}(t;x,y)\psi_2(y) \ \ .$$

\bigskip

\noi
A standard computation shows that

\begin{equation}\label{e:p2}
{\Ee}(t;x,y) = Q(t;x,y) + ({\Ee}*\Cc)(t;x,y),
\end{equation}

\bigskip

\noi where ${\Ee}*\Cc$ is a convolution given by

$$({\Ee}*\Cc)(t;x,y) = \int_0^tds\int_{M_R}dz \
{\Ee}(s;x,z)\Cc(t-s;z,y) ,$$

\bigskip

\noi
and the correction term $\Cc(t;x,y)$ is given by the formula

$$
\Cc(t;x,y) =
-{\frac{{\partial}^2\phi_1}{{\partial}u^2}}(x){\Ee}(t;x,y)\psi_1(y)$$
$$\qquad \qquad\ \ - {\frac{{\partial}\phi_1}{{\partial}u}}(x)
{\frac{{\partial}{\Ee}}{{\partial}u}}(t;x,y)\psi_1(y)
$$
$$
\qquad\qquad\ \ -{\frac{{\partial}^2\phi_2}{{\partial}u^2}}(x){
\Ee}(t;x,y)\psi_2(y) $$
 $$
 \qquad\qquad\ \ -
{\frac{{\partial}\phi_2}{{\partial}u}}(x)
{\frac{{\partial}{\Ee}}{{\partial}u}}(t;x,y)\psi_2(y).
$$

\bigskip
The choice of cut-off functions implies the following result:

\begin{lem}\label{l:err}
The ``error" term $\Cc(t;x,y)$ vanishes outside the cylinder $[\frac{1}{7} ,
\frac{6}{7}] \times Y$ and is equal to $0$ for
$d(x,y) < \frac{1}{7}$. Therefore, there exist positive constants $c_1,
c_2$ such that

\begin{equation}\label{e:err1}
\|\Cc(t;x,y)\| \le c_1e^{-c_3\frac{d^2(x,y)}{t}}.
\end{equation}
\end{lem}

\bigskip

We define the series

\begin{equation}\label{e:err2}
Q(t;x,y) + \sum_{n=1}^{\infty}(Q*\Cc_n)(t;x,y),
\end{equation}

\noi
where
$$\Cc_1 = C \ \ {\rm and}  \ \ \Cc_{n+1}(t;x,y) = \Cc_n*C \ \ .$$

\noi The elementary estimate

\begin{equation}\label{e:est0}
\|\Cc_n(t;x,y)\| \le \frac{c_1 {\rm vol}(y)t^{n-1}}{(n-1)!}
e^{-c_2\frac{d^2(x,y)}{t}}
\end{equation}

\bigskip

\noi
implies the absolute convergence of (\ref{e:err2}) and now the equality

$$\Ee(t;x,y) = Q(t;x,y) + \sum_{n=1}^{\infty}(Q*\Cc_n)(t;x,y)$$

\noi
is obvious. Proposition \ref{p:est}, jointly with
(\ref{e:err1}) and (\ref{e:est0})
gives us the following estimates on the kernels of the heat operators:

\bigskip

\begin{prop}\label{p:est1}
Let $\Ee(t;x,y)$ denote the kernel of one of the operators
$e^{-t\Dd_{\Pi_>}^2}$, $e^{-t\Delta_d}$,
$e^{-t\Delta_n}$, $e^{-t\Delta_{\pm}}$ on a manifold $M$ and let us denote
by $\Ff(t;x,y)$ the kernel
$\Dd\Ee(t;x,y)$. Assume that the corresponding Laplacian is an invertible
operator. Then there exist positive constants $c_1 , c_2, c_3$ such that

\begin{equation}\label{e:est1}
\qquad \|\Ee(t;x,y)\| \le
c_1t^{-\frac{n}{2}}e^{c_2t}e^{-c_3\frac{d^2(x,y)}{t}},
\end{equation}

$\ \ \|\Ff(t;x,y)\| \le
c_1t^{-\frac{n+1}{2}}e^{c_2t}e^{-c_3\frac{d^2(x,y)}{t}} \ .$
\end{prop}

\bigskip

\section{Duhamel's principle and the Adiabatic Limit}\label{s:str}

\bigskip
Now we want to analyze the behavior of the heat kernels in the
adiabatic limit. We start with the manifold $M$ with collar
neighborhood $N = [0,1] \times Y$ and we replace $M$ by $M_R$,
which is $M$ with $N$ replaced by $N_R = [0,R] \times Y$, a collar
of length $R$. To study the behavior of the heat kernels on $M_R$
we need the uniform estimates corresponding to the one we have
given in Proposition \ref{p:est}. To get them we use Duhamel's
Principle as in the previous Section, but now we take the
parameter $R$ into account.

More precisely, first we get the heat kernel ${\tilde
\Ee}_R(t;x,y)$ of the operator ${\tilde \Dd_R}^2$ on a manifold
${\tilde M_R}$. We obtain this kernel by gluing together the heat
kernel of $\Dd^2$ on the cylinder $(-\infty, +\infty) \times Y$
(restricted to $[-R,R] \times Y$) to the two copies of the heat
kernel  of $\Dd^2$ on $M$ (one for each end). The method described
in Section \ref{s:cyl} works in this case and the resulting kernel
${\tilde \Ee}_R(t;x,y)$ satisfies the estimate (\ref{e:est1}).

Now, we paste kernels together, but this time we make our
parametrix dependent on $R$. We
use the function $\rho(a,b)$ to define

$$
\phi_1 = 1 - \rho \left(\frac{5}{7}R , \frac{6}{7}R \right) , \; \psi_1 = 1
- \psi_{2},
$$
and
$$
\phi_2 = \rho \left(\frac{1}{7}R , \frac{2}{7}R \right) , \; \psi_2 =
\rho \left(\frac{3}{7}R , \frac{4}{7}R \right) ,
$$

\bigskip

\noi and introduce the corresponding functions on a manifold
$M_R$. We define $Q_R(t;x,y)$ a {\it ``parametrix''} for the heat
kernel ${\Ee}_R(t;x,y)$ (where again ${\Ee}_{cyl}(t;x,y)$ denotes
the heat kernel of one of our boundary problems)

\begin{equation}\label{e:r1}
Q_R(t;x,y) = \phi_1(x){\Ee}_{cyl}(t;x,y)\psi_1(y)
\end{equation}
$$
+ \phi_2(x){\tilde \Ee}_R(t;x,y)\psi_2(y) \ .$$

\bigskip
\noi
Again, we have

\begin{equation}\label{e:r2}
{\Ee}_R(t;x,y) = Q_R(t;x,y) + ({\Ee}_R*\Cc_R)(t;x,y),
\end{equation}

\bigskip
\noi where ${\Ee}_{R}*\Cc_R$ is a convolution and the correction
term $\Cc_{R}(t;x,y)$ is given by the formula from the previous
Section. The only difference is that cut-off functions depends on
$R$ . The crucial result is

\bigskip

\begin{thm}\label{l:r1}
The error term $\Cc_{R}(t;x,y)$ is equal to $0$ outside the
cylinder $[\frac{1}{7}R , \frac{6}{7}R] \times Y$. Moreover, it
is equal to $0$ if the distance between $x$ and $y$ is smaller
than $\frac{R}{7}$. As a result, there exist positive constants
$c_1, c_2, c_3$, such that the following estimate holds:

\begin{equation}\label{e:rr1}
\|{\Ee}_R(t;x,y)\| \le
c_1t^{-\frac{n}{2}}e^{c_2t}e^{-c_3{{\frac{d^2(x,y)}{t}}}}
\ \, .
\end{equation}
Moreover, the error term satisfies the estimate

\begin{equation}\label{e:ra1}
\|({\Ee}_R*\Cc_R)(t;x,x)\| \le
c_1e^{c_2t}e^{-c_3{{\frac{R^2}{t}}}}.
\end{equation}

\end{thm}
\bigskip
The proof goes exactly as before. We only sketch the proof of
(\ref{e:ra1}). In the following we use the vanishing of
$\Cc_R(t-s;z,x)$ for $d(x,z) > \frac{R}{7}$:
$$
\left \|({\Ee}_R*\Cc_R)(t;x,x) \right \|
$$
$$
= \left \|\int_0^tds\int_{M_R}{\Ee}_R(s;x,z)\Cc_R(t-s;z,x)dz \right \|
$$
$$
= \left \|\int_0^tds\int_{[\frac{R}{7} , \frac{6R}{7}] \times
Y}{\Ee}_R(s;x,z)\Cc_R(t-s;z,x)dz \right \|
$$
$$
\le  c_1e^{c_2t}\int_0^tds\int_{[\frac{R}{7} , \frac{6R}{7}]
\times Y} e^{-c_3{{\frac{td^2(x,z)}{s(t-s)}}}}dz
$$
$$
\le  c_1e^{c_2t}\int_0^tds\int_{[\frac{R}{7} , \frac{6R}{7}]
\times Y}e^{-c_4{{\frac{R^2}{t}}}}dz \qquad\qquad\qquad
$$
$$
\le c_1e^{c_2t}R{\cdot}{\rm vol}(Y)e^{-c_4{{\frac{R^2}{t}}}}\int_0^tds <
c_5e^{c_2t}e^{-c_6{{\frac{R^2}{t}}}} .
$$
\bigskip

Now we are able to show that the error contribution to the
$\z$-determinant for the {\it ``small''} time interval, meaning
$[0,R^{1-\e}]$, disappears in the adiabatic limit, i.e
\bigskip
\begin{cor}\label{c:rr1}
The following equality holds for small $\e >0$:

\begin{equation}\label{e:rr2}
\lim_{R \to \infty}\int_0^{R^{1-\e}}{\frac{dt}{t}}\int_{M_R}{\rm tr} \
({ \Ee}_R*\Cc_R)(t;x,x)dx = 0.
\end{equation}

\end{cor}

\bigskip

\begin{proof}
The result is an immediate consequence of Theorem \ref{l:r1}, because
$$
\left |\int_0^{R^{1-\e}}{\frac{dt}{t}}\int_{M_R}{\rm tr}
({\Ee}_R*\Cc_R)(t;x,x)dx \right |
$$
$$
\le
c_1\int_0^{R^{1-\e}}{\frac{e^{c_2t}}{t}}dt\int_{M_R}e^{-c_3{{\frac{R^2}{t}}}
}dx
$$
$$
\le c_1e^{c_2R^{1-\e}}e^{-c_3\over{2}R^{1 +
\e}}\int_0^{R^{1-\e}}{\frac{e^{-c_3{{\frac{R^2}{2t}}}}}{t}}dt\int_{M_R}dx
$$
$$
\le c_4Re^{c_2R^{1-\e}}e^{-c_3\over{2}R^{1 +
\e}}\int_0^{R^{1-\e}}{\frac{e^{-c_3{{\frac{R^2}{2t}}}}}{t}}dt
$$
$$
\le
c_5e^{-c_6R^{\e}} ,
$$
\noi
and (\ref{e:rr2}) follows easily.

\end{proof}
\bigskip
The meaning of the result is that as we take the adiabatic limit
the error contribution to the determinant can be neglected and we
are left only with the interior contribution and the cylinder
contribution. This, however holds only for a small time interval.
We will show in the next Section that the large time contribution
coming from the time interval $[R^{1-\e}, +\infty)$ can be
neglected.

\bigskip

\section{The small eigenvalues and the large time contribution}\label{s:eigen}

\bigskip

In this Section we explain why, in the adiabatic limit, we can
forget the contribution coming from the large time interval. Once
again we discuss only the simplest possible situation in which we
do not have to deal with small eigenvalues. We make the assumption
that the tangential operator $B$ is invertible. This condition
implies that there exists a constant $b > 0$ such that we have
only finitely many eigenvalues in the interval $[-b,b]$ for $R$
sufficiently large. To simplify further, in this exposition we
make one more assumption. We introduce the manifold $M_{\infty} =
((-\infty, 0] \times Y)\cup M$. The bundle $S$ and operator $\Dd$
extend naturally to $M_{\infty}$ and we assume that
\begin{equation}\label{e:1}
{\rm ker}_{L^2}\Dd = \{s \in L^2(M_{\infty};S) ; {\Dd}s = 0\} =
\{0\} \ \, .
\end{equation}

\noi
Assumption (\ref{e:1}) greatly simplifies the analysis of
the {\it Adiabatic Decomposition} of the $\z$-determinant. The reason is
that the operator $\Dd$ on $M_{\infty}$
has a unique closed, self-adjoint extension, which we denote by
$\Dd_{\infty}$ . This is a Fredholm operator (see
Section 6 of \cite{DoWo91}) and (\ref{e:1}) implies that the kernel of
$\Dd_{\infty}$ is equal to $\{0\}$. This implies the existence of
a positive constant $b$ such that for any spinor $s$ on $M_{\infty}$ we have

\begin{equation}
(\Dd^2s;s) \ge b\|s\|^{2}.
\end{equation}

\noi Let $\Delta_{R,\pm}$ denote the operator $\Delta_{\pm}$ on
the manifold $M_R$. We introduce similar notation for the other
boundary conditions. We also consider the operator ${\tilde
\Dd_R}$ , the Dirac operator on ${\tilde M_R}$ the double of
$M_R$. The operator ${\tilde \Dd_R}$ is the natural double of
$\Dd_R$ which is the Dirac operator $\Dd$ extended to $M_R$ .

\bigskip

\begin{prop}\label{p:1}
Let us assume that (\ref{e:1}) holds. Then there exists $R_0$ such that
$$
\mu > \frac{b}{2}
$$
\noi
for any eigenvalue $\mu$ of the operator $\Delta_{R,\pm}$, $\Delta_{R,d}$,
$\Delta_{R,n}$, $\Dd_{R,\Pi_>}^2$
or ${\tilde \Dd_R}^2$ and for any $R > R_0$.
\end{prop}

\bigskip

We do not present the proof of this result. It is not difficult but long
and technical.
The idea behind the result is easy to understand, however. Let $\la_0$
denote the smallest eigenvalue
of the operator $B^2$ and $\mu = \mu(R) < \la_0$ denote an eigenvalue of
one of the aforementioned Laplacians,
with the corresponding eigensection $\phi$. Assume that $\|\phi\| = 1$. We
can extend $\phi$ to the spinor
$\phi_{\infty}$ on $M_{\infty}$, which belongs to the domain of
$\Dd_{\infty}^2$. Moreover we can choose $\phi_{\infty}$
in a such a way that the $L^{2}$-norm of $\phi_{\infty}$ restricted to the
cylinder $M_{\infty} \setminus M_R$ is bounded as follows:

\begin{equation}\label{e:b/2}
\|\phi_{\infty}|_{M_{\infty} \setminus M_R}\|_{L^2}^2 \le
c_1e^{-c_2R} ,
\end{equation}

\noi for suitable positive constants $c_1, c_2$. Now the statement
of Proposition \ref{p:1} is an obvious consequence of min-max
principle. We refer to \cite{DoWo91} (Theorem 6.1) and
\cite{KPW94} (Proposition 2.1). A more general result was
published in \cite{Mu94}, Proposition 8.14.

All this implies that, in our {\it ``simple"} case, we can ignore the
large time contribution in the adiabatic limit.

\bigskip

\begin{prop}\label{c:lb1}
Let us assume (\ref{e:1})\,, then for any $\e >0$ the following
equality holds:
\begin{equation}\label{e:ltc}
\lim_{R \to \infty} \int_{R^{\e}}^{\infty}{\frac 1{t}}{\cdot}{\rm Tr} \
e^{-t\Delta_R}dt = 0
\, .
\end{equation}
\end{prop}

\bigskip

\begin{proof}
Assume that $R > R_0$ and let $\{\mu_k\}_{k=1}^{\infty}$ denote the set of
eigenvalues of $\Delta_R$\,. We have

$$
\int_{R^{\e}}^{\infty}{\frac 1{t}}{\cdot}{\rm Tr} \ e^{-t\Delta_R}dt =
\int_{R^{\e}}^{\infty}{\frac 1{t}}{\cdot}
\sum_{k = 1}^{\infty}e^{-t\mu_k}dt
$$
$$
=\int_{R^{\e}}^{\infty}{\frac 1{t}}{\cdot}\sum_{k =
1}^{\infty}e^{-(t-1)\mu_k}e^{-\mu_k}dt
$$
$$
< \int_{R^{\e}}^{\infty}{\frac 1{t}}e^{-(t-1){\frac{b}{2}}}{\cdot}{\rm Tr} \
e^{-\Delta_R}dt \ \ ,
$$

\bigskip

\noi where $b$ is the constant from Proposition \ref{p:1}. We now
have

$$\int_{R^{\e}}^{\infty}{\frac 1{t}}{\cdot}{\rm Tr} \ e^{-t\Delta_R}dt
$$
$$
<\int_{R^{\e}}^{\infty}{\frac 1{t}}e^{-(t-1){\frac{b}{2}}}{\cdot}{\rm Tr} \
e^{-\Delta_R}dt <
c_6R^{1-\e}{\cdot}e^{-c_7R^{\e}},
$$

\bigskip

\noi and the Proposition follows easily.
\end{proof}

\bigskip

\section{The decomposition of the $\z$-determinant of the Dirac
Laplacian}\label{s:lapl}

\bigskip
At last we are ready to discuss the decomposition of the $\z$-determinant.
The manifold $M$ is now an odd-dimensional closed manifold and we assume
that it has a decomposition $M_1 \cup M_2$\,, where $M_1$ and $M_2$
are compact manifolds with boundary such that

\begin{equation}\label{e:dec1}
M = M_1 \cup M_2  \ ,  \ M_1 \cap M_2 = Y = {\partial}M_1 =
{\partial}M_2.
\end{equation}
In this set-up $N$ denotes $N = [-1,1] \times Y$,  the bicollar
neighborhood  of $Y$ in $M$, and $N_R = [-R,R] \times Y$ is the
corresponding neighborhood of $Y$ in $M_R$. We denote by $\Dd$ the
Dirac operator on $M$ and $\Dd_i = \Dd|_{M_i}$. We want to find a
formula for the quotient
\begin{equation}
\label{e:d1}
\ \ \ \ \ \ \ \ \ \ \ \ \ \ \ \ \
\frac{{\rm det}_{\z}\Dd_R^2}
{{\rm det}_{\z}\Delta_{1,R,d}{\cdot}\Delta_{2,R,d}}
\ \,
\end{equation}

\noi or alternatively the difference

$$
{\rm ln} \ {\rm det}_{\z}\Dd_R^2
- {\rm ln} \ {\rm det}_{\z}\Delta_{1,R,d} - {\rm ln} \
{\rm det}_{\z}\Delta_{1,R,d}
$$
$$
= - \int_0^{\infty}\frac{dt}{t}\{{\rm Tr} \ e^{-t\Dd_R^2} - {\rm Tr} \
e^{-t\Delta_{1,R,d}} - {\rm Tr} \ e^{-t\Delta_{2,R,d}}\} .
$$

\bigskip

\noi It follows from the analysis presented in the previous
Sections that as $R \to \infty$ we can neglect the error terms and
study only the cylinder contribution to the difference ${\rm Tr} \
e^{-t\Dd_R^2} - {\rm Tr} \ e^{-t\Delta_{1,R,d}} - {\rm Tr} \
e^{-t\Delta_{2,R,d}}$. Moreover, on the cylinder the heat kernel
of the operator $\Dd_R^2$ can be replaced by the heat kernel
determined by the operator $-\partial_u^2 + B^2$. Modulo terms
which disappear as $R \to \infty$ we now have the equality
$$
{\rm Tr} \ e^{-t\Dd_R^2} - {\rm Tr} \ e^{-t\Delta_{1,R,d}}
- {\rm Tr} \ e^{-t\Delta_{2,R,d}}
$$
$$
= \int_{-R}^R\frac {1}{\sqrt{4{\pi}t}}{\rm Tr} \ e^{-tB^2}du
$$
$$
- \int_0^R\frac {1}{\sqrt{4{\pi}t}}(1-e^{-\frac {u^2}{t}}){\rm Tr} \
e^{-tB^2}du
$$
$$
- \int_{-R}^0\frac {1}{\sqrt{4{\pi}t}}(1-e^{-\frac {u^2}{t}}){\rm Tr} \
e^{-tB^2}du
$$
$$
= \frac {2}{\sqrt{4{\pi}t}}\int_0^Re^{-\frac {u^2}{t}}du{\cdot}{\rm Tr} \
e^{-tB^2}
$$
$$
= {\rm Tr} \ e^{-tB^2}{\cdot}\frac 1{\sqrt{\pi}}{\cdot}\int_0^{\frac R
{\sqrt{t}}}e^{-v^2}dv .
$$
\noi We obtain $\frac 1{2} {\rm Tr} \ e^{-tB^2}$ as $R \to \infty$ and
modulo minor technicalities we have proved

\begin{equation}\label{e:d2}
\lim_{R \to \infty}\{{\rm ln} \ {\rm det}_{\z}\Dd_R^2
- {\rm ln} \ {\rm det}_{\z}\Delta_{1,R,d}
\end{equation}
$\qquad\qquad - {\rm ln} \ {\rm det}_{\z}\Delta_{1,R,d}\} = \frac
1{2}{\cdot}{\rm ln} \ {\rm det}_{\z}B^2 ,$

\bigskip

\noi which yields our first adiabatic decomposition result

\bigskip

\begin{thm}\label{t:d1}
The following equality holds under the assumptions we have made:
\begin{equation}\label{e:d3}
\lim_{R \to
\infty}\frac{{\rm det}_{\z}\Dd_R^2}
{{\rm det}_{\z}\Delta_{1,R,d}{\cdot}\Delta_{2,R,d}}
={\sqrt{{\rm det}_{\z}B^2}}.
\end{equation}
\end{thm}

\vskip 1cm

We work out the case of the Neumann condition in the same way. The only
difference is the sign of the contribution and
therefore we obtain

\begin{equation}\label{e:n1}
\lim_{R \to
\infty}\frac{{\rm det}_{\z}\Dd_R^2}
{{\rm det}_{\z}\Delta_{1,R,n}{\cdot}\Delta_{2,R,n}}
=\frac 1{ \sqrt{{\rm det}_{\z}B^2}}.
\end{equation}

\bigskip

\noi This method also works in the case of the chiral boundary
condition. In this case the Neumann contribution cancels out the
Dirichlet contribution and as the result we have the formula

\begin{equation}\label{e:chsplit1}
\lim_{R \to
\infty}\frac{{\rm det}_{\z}\Dd_R^2}
{{\rm det}_{\z}\Delta_{1,R,+}{\cdot}\Delta_{2,R,+}}
= 1.
\end{equation}

\bigskip

\noi This formula was somehow the first we noticed and we used it
to obtain the corresponding result for the Atiyah--Patodi--Singer
condition (see \cite{JKPW2})
\bigskip

\begin{thm}\label{t:chiralsplit2}
The following equality holds in the case of Atiyah--Patodi--Singer condition:
\begin{equation}\label{e:chiralsplit3}
\lim_{R \to \infty}
{\frac{{\rm det}_{\z}\Dd_R^2}{{\rm det}_{\z}\Dd_{1,R,\Pi_<}^2{\cdot}
{\rm det}_{\z}\Dd_{2,R,\Pi
_>}^2}} = 2^{-\z_{B^2}(0)}.
\end{equation}
\end{thm}
\bigskip

We refer to \cite{JKPW2} for the details of the proof. This ends the
discussion of the decomposition of the $\z$-determinant of
the Dirac Laplacian in case we do not have to deal with the small eigenvalues.

\bigskip

\section{The splitting of the $\eta$-invariant}\label{s:eta}

\bigskip
Here we discuss the decomposition of the phase of the $\z$-determinant. Let
us first observe that the $\eta$-function of the
Atiyah--Patodi--Singer boundary problem shares the properties of the
$\eta$-function of the Dirac operator on a closed manifold.
This is due to the fact that the boundary does not create any new
singularities of the $\eta$-function. The singularities are created
by the small
time asymptotics of the trace ${\rm Tr} \ \Dd e^{-t\Dd^2}$.
Once again using Duhamel's principle we see that we have to study the trace
of the Heat Kernel on the cylinder. The kernel has the form
$$
G(\partial_u + B)\Ee_>(t;(u,x),(v,y)) \ \ ,$$

\bigskip

\noi where $\Ee_>(t;(u,x),(v,y))$ is the kernel of the operator on
the cylinder (see the formulas (\ref{e:g1})). We only need to
notice that it has the form

$$\sum g_n(t;u,v)\phi_n(x)\tensor \phi_n(y) \ \ ,$$

\bigskip
\noi
where $\{\phi_n\}$ is the orthonormal basis of eigenspinors of $B$. Let us
recall that we can choose this basis
in such a way that $G\phi_n = \phi_{-n}$. Now it follows from (\ref{e:GB})
that the trace of

$$G(\partial_u + B)\Ee_>(t;(u,x),(v,y))$$

\noi in the $Y$-direction is equal to $0$.

\bigskip
It is not difficult to see that in fact not only the $\eta$-invariant, but the
$\z$-determinant is well-defined on the whole Grassmannian
$Gr_{\infty}^*(\Dd)$ . This happens because the
$\eta$- and $\z$-functions behave nicely on this
particular space of boundary conditions.
We start with a more precise description of the {\it Unitary
Twist}, which we already encountered in
Section 3.

\bigskip

\begin{lem}\label{l:path} For any $P \in Gr^*_{\infty}(\Dd)$ there exists a
smooth path
$\{g_u\}_{0 \le u \le 1}$ of unitary operators on $L^2(Y;S|Y)$ which
satisfies
\[ Gg_u = g_uG \ \ and \ \ g_u - Id \ \  has \ a \ smooth \ kernel,
\]

\noi
such that $g_1 = Id$ and the path  $\{P_u = g_u{\Pi_>}g_u^{-1}\} \subset
Gr^*_{\infty}(\Dd)$
connects $P_0 = P$ with $P_1 = \Pi_>$.
\end{lem}

\bigskip
We can always assume that the path $\{g_u\}$ is constant on subintervals
$[0,1/4]$ and $[3/4,1]$.
We introduce $U$, a unitary operator on $L^2(M;S)$. The operator $U$ is
equal to the $Id$ on the complement of the collar $N$ and
$$
U|_{\{u\}\times Y} = g_u .
$$

\bigskip

 The following Lemma introduces the {\it Unitary Twist}, which
allows us to replace the operator $\Dd_P$ by a modified operator
$\Dd + \Rr$ subject to the  boundary condition $\Pi_{>}$. This
makes an explicit construction of the heat kernels on a cylinder
possible.
\bigskip
\begin{lem}\label{l:unitary}
The operators $\Dd_P$ and $\Dd_{U,\Pi_>} = (U^{-1}{\Dd}U)_{\Pi_>}$ are
unitarily equivalent operators.
\end{lem}

\begin{proof}
Let $\{f_k;\m_k\}_{k \in \Z}$ denote a spectral resolution of the operator
$\Dd_P$\,. This
means that for each $k$ we have

$$
\Dd{f}_k = \m_kf_k \ \ and \ \ P(f_k|Y) = 0 \ \ .
$$
\noi
This implies
$$
U^{-1}{\Dd}U(U^{-1}f_k) = \m_k(U^{-1}f_k)$$
\noi
and
$$
\Pi_{\s}((U^{-1}f_k)|Y) = g_0^{-1}P(f_k|Y) = 0 \ \ .$$
\bigskip
\noi
hence $\{U^{-1}f_k;\m_k\}$ is a spectral resolution of
$(U^{-1}{\Dd}U)_{\Pi_>}$ .

\end{proof}

\bigskip

In the collar $N$ , we have the formulas

\[
\ \ \ \
U^{-1}{\Dd}U = \Dd +
GU^{-1}{{{\partial}U}\over{{\partial}u}} + GU^{-1}[B,U] \ , \]

\noi
and
\[
\ \ \ \ \ \ \ \ \ \
U^{-1}{\Dd}^2U = \Dd^2 -
2U^{-1}{{{\partial}U}\over{{\partial}u}}\partial_u \]
\[
\ \ \ \ \ \ \ \ \ \
 -U^{-1}{{{\partial}^2U}\over{{\partial}u^2}} + U^{-1}[B^2,U] \ ,
\]
\bigskip

\noi
which, restricted to the collar $[0,1/4] \times Y$, give

\begin{equation}\label{e:udu}
\ \ \ \ \ \ \ \ \ \ \ \ \
U^{-1}{\Dd}U = \Dd + GU^{-1}[B,U], \ \,
\end{equation}

and
$$
U^{-1}{\Dd}^2U = {\Dd}^2 + U^{-1}[B^2,U] \ .
$$

\bigskip

\noi
It follows from Lemma \ref{l:unitary} that we can study the operator
$\Dd_{U,\Pi_>}$ instead of
the operator $\Dd_P$. Again, it is enough to study the small time
asymptotics of the trace of the Heat operator on the cylinder.
This all tells you that up to an exponentially small error we have to study
the trace of the operator
kernel of the operator

\begin{equation}\label{e:cyl1}
(G(\partial_u + B) + \Kk_1)e^{-t(-\partial_u^2 + B^2 + \Kk_2)_{\Pi_>}}
\ \,
\end{equation}

\noi
where

\[
\ \ \ \ \
\Kk_1 = GU^{-1}[B,U] \ \ and \ \ \Kk_2 = U^{-1}[B^2,U] \
\ ,
\]

\bigskip
\noi
to study meromorphic extension of the $\eta$-function, and simply the trace

$$Tr \ e^{-t(-\partial_u^2 + B^2 + \Kk_2)_{\Pi_>}}$$

\noi
to learn about $\z$-function. Let us observe that $\Kk_1$
anticommutes and $\Kk_2$ commutes with
the involution $G$. The point is that $\Kk_1$ and $\Kk_2$ are smoothing in
the $Y$-direction. Hence when we study the trace
${\rm Tr} \ e^{-t\Dd_{U,\Pi_>}^2}$ we can easily show that
\begin{equation}\label{e:zz}
|{\rm Tr} \ e^{-t\Dd_{U,\Pi_>}^2} - {\rm Tr} \ e^{-t\Dd_{\Pi_>}^2}| <
c\sqrt{t}.
\end{equation}
\bigskip
Modulo exponentially small term this difference is equal to the sum given
by Duhamel's principle. The first term here is
$$
{\rm Tr}(\Ee_{\Pi_>}*\Kk_2\Ee_{\Pi_>})(t) =
$$
$$
\int_0^t{\rm Tr}\Ee_{\Pi_>}(s)\Kk_2\Ee_{\Pi_>}(t-s)ds
= t{\cdot}{\rm Tr} \
\Kk_2\Ee_{\Pi_>}(t) .
$$
\bigskip
\noi
The operator $\Kk_2$ smoothes things out in the $Y$-direction so the only
singularity left is in the normal direction and we have
$$
|{\rm Tr} \ \Kk_2\Ee_{\Pi_>}(t)| < \frac c{\sqrt{t}} \ \ ,$$
\bigskip
\noi and now (\ref{e:zz}) follows. Details are given in
\cite{KPW99}. Straightforward computations show that

$$\z_{\Dd_{\Pi_>}^2}(0) = 0$$

\bigskip
\noi
(see \cite{YLKW01} and \cite{JKPW2}). This fact combined
with the estimate (\ref{e:zz}) gives the following result:

\bigskip

\begin{prop}\label{p:zz}
(\cite{YLKW01}).

$$\z_{\Dd_P^2}(0) = 0$$

\bigskip
\noi
for any $P \in Gr_{\infty}^*(\Dd)$, such that $\Dd_P$ is invertible.
\end{prop}

\bigskip

\begin{proof}
We have
$$
\z_{\Dd_P^2}(0) = \z_{\Dd_P^2}(0) - \z_{\Dd_{\Pi_>}^2}(0)
$$
$$
=\lim_{s \to 0}\frac 1{\G(s)} \int_0^{\infty}t^{s-1}{\rm Tr}(e^{-t\Dd_P^2} -
e^{-t\Dd_{\Pi_>}^2})dt
$$
$$
= \lim_{s \to 0} s\int_0^1t^{s-1}{\rm Tr}(e^{-t\Dd_P^2}
- e^{-t\Dd_{\Pi_>}^2})dt.
$$
Now we use (\ref{e:zz}) to end the proof, i.e.
$$
|\z_{\Dd_P^2}(0) - \z_{\Dd_{\Pi_>}^2}(0)| \le c\lim_{s \to 0}
s\int_0^1t^{s- \frac 1{2}}dt = 0 .
$$
\end{proof}

\bigskip

\noi Similarly we show
\begin{equation}\label{e:eta1}
|{\rm Tr} \Dd_Pe^{-t\Dd_P^2} - {\rm Tr}
\Dd_{\Pi_>}e^{-t\Dd_{\Pi_>}^2}| \le c,
\end{equation}
\bigskip
\noi
which implies the following result:
\bigskip

\begin{thm}\label{t:eta1}
For any $P \in Gr_{\infty}^*(\Dd)$ the function $\eta_{\Dd_P}$ is a
holomorphic function of $s$ for $Re(s) > -1$ .
In particular we have the equality
\begin{equation}\label{e:eta11}
\eta_{\Dd_P}(0) = \frac 1{\sqrt{\pi}}\int_0^{\infty}\frac
1{\sqrt{t}} {\rm Tr} \Dd_Pe^{-t\Dd_P^2}dt.
\end{equation}
\bigskip
\noi
Moreover, let $\{P_r\}$ denote a smooth family of projections from
$Gr_{\infty}^*(\Dd)$. The variation of the
$\eta$-invariant of the family $\{\Dd_{P_r}\}$ is given by the formula

\begin {equation}\label{e:eta12}
\qquad {d\over{dr}}\{\eta_{\Dd_{P_r}}(0)\}|_{r=0} =
\end{equation}
$\qquad\qquad -{2\over{\sqrt{\pi}}}{\cdot} \lim_{\e \to
0}\sqrt{\e}{\cdot}{\rm Tr} \ {\dot \Dd}_{P_0}e^{-{\e}\Dd_{P_0}^2}
$

\bigskip

\noi
where ${\dot \Dd}_0 =\frac d{dr}\{U_r^{-1}{\Dd}U_r\}|_{r=0}$ .
\end{thm}

\bigskip

It follows that the $\z$-determinant gives a well-defined smooth function

$$
P \to {\rm det}_{\z}\Dd_P$$

\bigskip
\noi
on the Grassmannian $Gr_{\infty}^*(\Dd)$ .

\bigskip
Let us discuss the decomposition of the $\eta$-invariant and its
dependence on the choice of the boundary conditions on $M_1$ and
$M_2$. If we fix Atiyah--Patodi--Singer conditions on both $M_1$
and $M_2$, then we know that there is no boundary contribution. We
repeat the analysis from the case of the Dirac Laplacian. This
gives

\bigskip

\begin{thm}\label{t:edc1}
The following formula holds under the assumptions we have made on the
operator $\Dd$:

\begin{equation}\label{e:edc1}
\lim_{R \to \infty}\{\eta_{{\Dd_R}}(0) - \eta_{\Dd_{1,R,\Pi_<}}(0)
- \eta_{\Dd_{2,R,\Pi_>}}(0)\} = 0.
\end{equation}
\end{thm}

\bigskip
Theorem \ref{t:edc1} corresponds to the results on the
decomposition of the modulus of the $\z$-determinant. However, the
$\eta$-invariant is a much more rigid invariant than ${\rm
det}_{\z}\Dd^2$. First of all the variation of the
$\eta$-invariant is given by the local formula. This gives an
immediate corollary. Namely (\ref{e:edc1}) holds independently of
$R$ . It is easy to understand what is going on. We replace small
part of the cylinder (away from the boundary) by a longer piece.
In this way we change the operator $\Dd_R$. Now the variation of
the $\eta$-invariant feels only what is going on at the given
point and from the local point of view the process corresponds to
the study of the operator $G(\partial_u + B)$ on the manifold
$S_R^1 \times Y$, where $S_R^1$ is the circle with radius $R$.
This operator has a symmetric spectrum. If $\phi(u,y)$ is an
eigenspinor corresponding to the eigenvalue $\la$, then

$${\bar \phi}(u,y) = G\phi(2{\pi}R-u,y)$$

\noi
is the eigensection corresponding to the eigenvalue $-\la$. Therefore for any
$R$ the $\eta$-function on $S_R^1$ disappears. This gives

\begin{equation}\label{e:edc2}
\eta_{{\Dd_R}}(0) - \eta_{\Dd_{1,R,\Pi_<}}(0) -
\eta_{\Dd_{2,R,\Pi_>}}(0) = 0.
\end{equation}

Now we want to relax the assumptions on the small eigenvalues. There is no
problem with the situation
in which the operators $\Dd_i$ have nontrivial $L^2$-kernel when extended
to $M_{i,\infty}$.
We simply modify the operator $\Dd_i$ by adding the orthogonal
projection  onto the $L^2$-kernel. The only thing which may vary
during this process
is that a finite number of eigenvalues might change the sign, or become
zero modes. For this reason we can only discuss the equalities
$mod \ \Z$ . The computation of the integer contribution
has to be done separately and uses different methods.
Similar modification lead to the relaxing of the condition of the
invertibility of the tangential operator $B$ .

Let us assume that the operator $B$ has non-trivial kernel. The involution $G$
(see (\ref{e:GB})) restricted to $ker(B)$ defines
a symplectic structure on this subspace of $L^2(Y;S|Y)$ and the
{\it Cobordism Theorem for Dirac
Operators} (see for instance \cite{BoWo93}, Corollary 21.16) implies

\[
\ \ \ \ \ \ \ \ \ \
{\rm dim} \ {\rm ker}(B^+) = {\rm dim} \ {\rm ker}(B^-)  \ \ . \]

\noi
This last equality shows the existence of Lagrangian subspaces of $ker(B)$.
We choose such
a subspace $W$ and let $\s : L^2(Y;S|Y) \to L^2(Y;S|Y)$ denote the
orthogonal projection of
$L^2(Y;S|Y)$ onto $W$. Let $\Pi_>$ denote the orthogonal projection of
$L^2(Y;S|Y)$ onto the
subspace spanned by eigenvectors of $B$ corresponding to the positive
eigenvalues. Then

\begin{equation}\label{e:pr}
\Pi_{\s} = \Pi_> + \s  \in Gr^*_{\infty}(\Dd) \ \,
\end{equation}

\noi
gives an element of $Gr_{\infty}^*(\Dd)$, which is a finite-dimensional
perturbation of the Atiyah--Patodi--Singer condition.

Let $\Pi_{\s}$ denote
a projection given by Formula (\ref{e:pr}). We repeat our analysis again
and obtain the following formula:

\begin{equation}\label{e:add2}
\eta_{\Dd}(0) = \eta_{{\Dd_{1_{Id - \Pi_{\s}}}}}(0) +
\eta_{{\Dd_{2_{\Pi_{\s}}}}}(0) \ \ {\rm mod} \ \Z \ \,.
\end{equation}

Now we introduce the formula for the variation of the
$\eta$-invariant under a change of boundary condition. The correct
approach to the computation was proposed by Lesch and
Wojciechowski (see \cite{LeWo96}).

Let $P \in Gr_{\infty}^*(\Dd)$ and let us choose a path
$\{P_r\}_{0 \le r \le 1} \subset Gr_{\infty}^*(\Dd)$ such that $P_0 =
\Pi_{\s}$ and $P_1 = P$.
There exists a smooth family $\{g_r\}$
of unitary
operators of the form \hskip 2mm $Id|_{(S|Y)} + smoothing \ operator$
which commutes with $G$
and such that

\[
\ \ \ \ \ \ \ \ \ \ \
g_0 = Id \ \ {\rm and} \ \ g_1\Pi_{\s}g_1^{-1} = P
\ \, .
\]

\bigskip

Next, we choose a smooth non-increasing function $\g(u)$ such that

\[
\g(u) = 1 \ \ {\rm for} \ u < 1/4, \ g(u) = 0 \ \ {\rm for} \ u > 3/4,
\]

\noi
and for each $0 \le r \le 1$ use the family

\begin{equation}\label{e:un}
\ \ \ \ \ \
g_{r,u} = g_{r\g(u)} \ \ {\rm for} \ \ 0 \le u \le 1 \ \, ,
\end{equation}

\bigskip

\noi
in order to construct a corresponding unitary operator $U_r$ on $M_2$ .
The operator $\Dd_{2_{U_r,\s}}$ is unitarily equivalent to the operator
$\Dd_{2_{P_r}}$. The variation of the
$\eta$-invariant is given by the standard formula (\ref{e:eta12}), which
allows us to prove the next result.

\bigskip

\begin{thm}\label{t;var}
For any $P \in Gr_{\infty}^*(\Dd)$, and any path $g = \{g_{r,u}\}$
connecting $\Pi_{\s}$ with $P$,
as described above, the following formula holds:
$$
\eta_{\Dd_{2_{P}}}(0) - \eta_{\Dd_{2_{\Pi_{\s}}}}(0)
$$
\begin{equation}\label{e:var2}
= -{1\over{\pi}}
\int_0^1dr\int_0^1du \
{\rm Tr} \
G \left ({\dot {g^{-1}{{{\partial}g}
\over{{\partial}u}}}}\right)|_{r} \ \ {\rm mod} \ \Z,
\end{equation}
\noi
where $({\dot {g^{-1}{{{\partial}g}\over{{\partial}u}}}})|_{r_0} =
{d\over{dr}}({g^{-1}{{{\partial}g}\over{{\partial}u}}})|_{r=r_0}.$
\end{thm}
\bigskip
\begin{proof}
We show that

\begin{equation}\label{e:var3}
{2\over{\sqrt{\pi}}}\lim_{\e \to 0} \sqrt{\e} {\rm Tr} \
(d(\Dd_{2_{U_r,\s}})/dr)|_{r=r_0}e^{-\e\Dd_{2_{U_{r_0},\s}}^2},
\end{equation}
$$
= {1\over{\pi}}\int_0^1 {\rm Tr} \ G \left({\dot
{g^{-1}{{{\partial}g}\over{{\partial}u}}}}\right)|_{r=r_0}du.
$$
We have
$$
({\dot {U^{-1}{\Dd}U}}) = G{\dot {U^{-1}{{{\partial}U}\over{{\partial}u}}}} +
G[U^{-1}BU , U^{-1}{\dot U}]$$
$$\qquad \qquad = G{\dot
{g^{-1}{{{\partial}g}\over{{\partial}u}}}} + G[g^{-1}Bg ,
g^{-1}{\dot g}].
$$
Thus
\[
\lim_{\e \to 0} \sqrt{\e}{\cdot}{\rm Tr} \ (d(\Dd_{2_{U_r,\s}})/dr)|_{r=r_0}
e^{-\e\Dd_{2_{U_{r_0},\s}}^2}
\]
\bigbreak
\noi
contains two terms. Let us start with
\[
\lim_{\e \to 0}
\sqrt{\e}{\cdot}{\rm Tr} \
G[g^{-1}Bg , g^{-1}{\dot g}]e^{-\e\Dd_{2_{U_{r_0},\s}}^2}\,.
\]
Once again we use {\it Duhamel's Principle} and replace
$e^{-\e\Dd_{2_{U_{r_0},\s}}^2}$ by the operator
$exp(-t(-\partial_u^2 + B^2 + \Kk_2)_{\s})$ on the cylinder. The
point here is that the kernel of this operator commutes with $G$
and the operator $[g^{-1}Bg , g^{-1}{\dot g}]$ anticommutes with
the involution $G$. It follows that

\[
Tr \ G[g^{-1}Bg , g^{-1}{\dot g}]e^{-\e\Dd_{U_r,\s}^2} = O(e^{-c/\e})
\ \,,
\]

\noi
and one is left with
\[
{2\over{\sqrt{\pi}}}\lim_{\e \to 0}\sqrt{\e}{\cdot}{\rm Tr} \
G \left({\dot
{U^{-1}{{{\partial}U}\over{{\partial}u}}}}\right)
e^{-\e\Dd_{2_{U_{r_0},\s}}^
2}
\ \, .
\]

\bigskip
\noi The term $G({\dot
{U^{-1}{{{\partial}U}\over{{\partial}u}}}})|_{r=r_0} = G({\dot
{g^{-1}{{{\partial}g}\over{{\partial}u}}}})|_{r=r_0}$ is supported
in $[1/4,3/4] \times Y$, and so we replace the kernel of the
operator $e^{-\e\Dd_{2_{U_{r_0},\s}}^2}$ by the kernel of the
operator $exp(-\e(-\partial_u^2 + B^2))$ on the infinite cylinder
$(-\infty,+\infty) \times Y$. Now we have

\[
\quad {2\sqrt{\e}\over{\sqrt{\pi}}}{\cdot}{\rm Tr} \ G \left({\dot
{g^{-1}{{{\partial}g}\over{{\partial}u}}}}\right)|_{r_0}
e^{-\e\Dd_{2_{U_r,\s}}^2}
\]
\[
={2\sqrt{\e}\over{\sqrt{\pi}}} \int_0^1du \ {\rm Tr}_Y \ G
\left({\dot
{g^{-1}{{{\partial}g}\over{{\partial}u}}}}\right)|_{r_0}
g^{-1}e^{-\e(-\partial_u^2 + B^2)}g
\]
\[
={2\sqrt{\e}\over{\sqrt{\pi}}} \int_0^1du \ {\rm Tr}_Y \ G
\left({\dot
{g^{-1}{{{\partial}g}\over{{\partial}u}}}}\right)|_{r_0}
e^{-\e(-\partial_u^2 + B^2)}
\]
\[
={2\sqrt{\e}\over{\sqrt{\pi}}}{1\over{\sqrt{4\pi\e}}} \int_0^1du \
{\rm Tr}_Y \ G \left({\dot
{g^{-1}{{{\partial}g}\over{{\partial}u}}}}\right)|_{r_0}e^{-{\e}B^2}
\]
\[
={1\over{\pi}}\int_0^1du {\rm Tr}_Y \ G \left({\dot
{g^{-1}{{{\partial}g}\over{{\partial}u}}}}\right)|_{r_0}
e^{-{\e}B^2}\]

\bigskip

\noi so that

\[
\quad {2\over{\sqrt{\pi}}}\lim_{\e \to 0}\sqrt{\e}{\cdot}{\rm Tr}
\ G \left({\dot
{g^{-1}{{{\partial}g}\over{{\partial}u}}}}\right)|_{r_0}
e^{-\e\Dd_{2_{U_r,\s}}^2}
\]
$$ ={1\over{\pi}}\int_0^1du \ {\rm Tr}_Y \ G \left({\dot
{g^{-1}{{{\partial}g}\over{{\partial}u}}}}\right)|_{r_0}.\qquad\qquad\qquad\qquad
$$
\end{proof}
\bigskip
\begin{rem}
If we assume that $g_r(u)$ is given by the formula
\[
\ \ \ \ \ \
g_r(u) = \pmatrix{
      Id & 0 \cr
0 & {\rm exp}(ir{\g}(u))\Theta \cr}
\ \, ,
\]

\bigskip

\noi
where $\Theta : C^{\infty}(Y;S^-|Y) \to C^{\infty}(Y;S^-|Y)$ is a self-adjoint
operator with a smooth kernel, then our formula has a very nice and simple form

$$\eta_{\Dd_{2_{P}}}(0) - \eta_{\Dd_{2_{\Pi_{\s}}}}(0)\qquad\qquad \qquad\qquad\qquad\qquad$$

\begin{equation}\label{e:etabc}
= -{1\over{\pi}}\int_0^1dr\int_0^1du \ {\g}'(u){\rm Tr} \ \Theta
= {{{\rm Tr}
\ \Theta}\over{\pi}}
\end{equation}
${\rm mod} \; \Z$.
\noi
This is the formula obtained by Lesch and Wojciechowski for the
finite-dimensional perturbation of the Atiyah--Patodi--Singer
condition (see \cite{LeWo96}).
\end{rem}

\bigskip

\begin{cor}\label{c:var1}
Let $P_1, P_2 \in Gr^*_{\infty}({\Dd})$ , then

$$\eta_{\Dd_{2_{P_1}}}(0) - \eta_{\Dd_{2_{P_2}}}(0)\qquad\qquad\qquad\qquad\qquad\qquad$$
\begin{equation}\label{e:var11}
= -{1\over{\pi}}
\int_0^1dr\int_0^1du \
{\rm Tr} \
G \left({\dot {g^{-1}{{{\partial}g}\over{{\partial}u}}}}\right)|_{r}
\ \ {\rm mod} \ \Z
\end{equation}

\noi
where $\{g_{r,u}\}$ is any family connecting $P_1$ with $P_2$ in the way
described above
(see (\ref{e:un})).
\end{cor}

\bigskip

\begin{cor}\label{c:ph}
The variation of the $\eta$-invariant
${d\over{dr}}(\eta_{\Dd_{2_{P_r}}}(0))|_{r=0}$ does not
depend on the choice of the base projection $P = P_0$. It depends only on
the family of
unitary operators $\{g_r\}$.
\end{cor}

\bigskip

This result plays an important role in the proof of equality
of the $\z$-determinant and the $\Cc$-determinant.

\bigskip

\begin{thm}\label{t:main1}
For any $P_1, P_2 \in Gr_{\infty}^*(\Dd)$ one has the following formula:

\begin{equation}\label{e:add1}
\eta_{\Dd}(0) =\eta_{{\Dd_{1_{Id - P_1}}}}(0) +
\eta_{{\Dd_{2_{P_2}}}}(0)
\end{equation}
$ \qquad\qquad\qquad\qquad\qquad + \eta(P_1,P_2)(0) \ \ {\rm mod}
\ \Z ,$

\bigskip

\noi where $\eta(P_1,P_2)(0)$ denotes the eta-invariant of the
operator $G(\partial_u + B)$ on $[0,1] \times Y$ subject to the
boundary condition equal to $P_1$ at $u=0$ and $Id-P_2$ at $u=1$.
\end{thm}

\bigskip

We need to explain the appearance of the middle term. We start with
the equality

$$\eta_{\Dd}(0) =\eta_{{\Dd_{1_{Id - \Pi_{\s}}}}}(0)+ \eta_{{\Dd_{2_{\Pi_{\s}}}}}(0)\qquad\qquad$$
$ \qquad\qquad\qquad\qquad\qquad+ \eta(\Pi_{\s}, \Pi_{\s})(0)  \ \
{\rm mod} \ \Z .
$

\bigskip

\noi The last term on the right side is equal to $0$ by virtue of
the natural symmetry described earlier in this Section. Now we
vary the boundary conditions replacing $Id - \Pi_{\s}$ by $Id-P_1$
on $M_1$ and $\Pi_{\s}$ by $P_1$ on the left end of the cylinder.
Then we replace $\Pi_{\s}$ by $P_2$ on $M_2$ and $Id-\Pi_{\s}$ by
$Id-P_2$ on the right end of the cylinder. The total variation of
the $\eta$-invariant under these changes is equal to $0$ (${\rm
mod} \ \Z$).

\bigskip

\section{Some remarks on the dependence on $\R$}\label{s:lapl}

\bigskip

In general, as one might expect, the ratios of the
$\z$-determinant discussed in this paper depend on the length of
the cylinders connecting two different parts of the manifolds. We
made explicit computations in which the Fredholm determinant shows
up and it is easy to see its explicit $R$-dependence (unpublished
work of the authors). However, here we  want to study the case in
which the ratio is $R$-independent. This situation brings up
another nice adiabatic picture to the story. The approach is based
on the work of L. Nicolaescu (see \cite{N95}).

\bigskip

We now denote by $M_R$ the manifold

$$M_{R} = M_2 \cup [-R,0] \times Y \ \ .$$

\bigskip

\noi
We have a $1$-parameter family of Cauchy data spaces of $\Dd_R$,
$\Lambda^{R}(D)$.
For any non-negative real number $\nu$, we define
$$
H_{\nu}={\rm span}_{L^2}\{\phi \ | \  B\phi=\lambda \phi
\ \ {\rm and} \ \
|\lambda|\le \nu\} \ \ ,
$$
$$
H^{\nu}_{<}={\rm span}_{L^2}\{\phi \ | \ B\phi=\lambda \phi \ \
{\rm and} \ \ \lambda
< \nu \} \ \ ,
$$
$$
H^{\nu}_{>}={\rm span}_{L^2}\{\phi \ | \ B\phi=\lambda \phi \ \
{\rm and} \ \ \lambda
> \nu \} \ \ .
$$
\bigskip

\noi It is well known that $\Lambda^0(D)$ and $H^{\nu}_{<}\oplus
U$ are the Fredholm pair for any $\nu\in\mathbb{R}$, and any
finite dimensional subspace $U \subset L^2(Y,S|_Y)$. Hence there
exists a number $\nu_0$ such that

$$\Lambda^0(D) \cap H^{\nu_0}_{<}=0.$$

\bigskip
\noi
The smallest such $\nu_0$ is called the non-resonance level of $D$. The
symplectic reduction
of $\Lambda^R(D)$ to $H_{\nu}$, which is defined by

$$L^R_{\nu}:=\frac{\Lambda^R(D)\cap(H_{\nu}\oplus H^{\nu}_{<})}
{H^{\nu}_{<}}\subset H_{\nu}$$

\bigskip
\noi
is the Lagrangian subspace of $H_{\nu}$.
Since the tangential operator $B$ preserves $H_{\nu}$, we can form the
$1$-parameter family of finite-dimensional operators

$$ e^{-R B_{\nu}}: H_{\nu} \quad \to \quad H_{\nu},
$$
\bigskip

\noi where $B_{\nu}$ is the restriction of $B$ to $H_{\nu}$. We
need the following description of the dynamics of the Cauchy data
space $\Lambda^R(D)$:

\bigskip

\begin{prop}\label{p:dyn-cd}
(\cite{N95}). For $\nu\ge\nu_0$, as $R\to\infty$,
$$
\Lambda^R(D) \ \to \ L^{\infty}_{\nu}\oplus H^{-\nu}_{<},
$$

\noi where

$$
L^{\infty}_{\nu}= \lim_{R\to\infty} L^R_{\nu} = \lim_{R\to\infty} e^{-R
B_{\nu}} L^R_{\nu}.
$$
\end{prop}

\bigskip

\noi
Nicolaescu's proposition leads to the following interesting result.

\bigskip

\begin{prop}\label{p:R-dep}
Given a couple of boundary conditions $(P_1,P_2)=(\Pi_{>}+
\sigma_1, \Pi_{>}+\sigma_2)$ where $\sigma_1$, $\sigma_2$ are the
orthogonal projections to the Lagrangian subspaces $L_1, L_2$ of
 $\ H_{0} = {\rm ker}\ B$, we assume that ${\rm
ker}(\Dd_R)_{P_1}={\rm ker}(\Dd_R)_{P_2}=0$. Then the quotient
$$
\frac{{\rm det}_{\z}(\Dd_{R})_{P_1}^2}
{{\rm det}_{\z}(\Dd_{R})_{P_2}^2}
$$
\noi
does not depend on $R$ .
\end{prop}

\bigskip

\begin{proof}
The proof of this proposition is an application of the following
Scott--Wojciechowski formula, Proposition 4.1 in \cite{SSKPW299},

\bigskip

$\frac{{\rm det}_{\z}(\Dd_{R})_{P_1}^2}
{{\rm det}_{\z}(\Dd_{R})_{P_2}^2}=$
\begin{equation}\label{e:sw}
\qquad |{\rm
det}_{Fr}U_{P_2}(U_{P_1})^{-1}\Ss_R(P_1)\Ss_R(P_2)^{-1}|^{2},
\end{equation}
\bigskip
\noi
where ${\rm det}_{Fr}$ is the Fredholm determinant, and

$$
U_{P_2}(U_{P_1})^{-1}:{\rm Range}(P_1)\to {\rm Range}(P_2)
$$

\bigskip
\noi is an unitary map which depends only on $P_1,P_2$. The
operators $U(P)$ and $S(P)$ were introduced in Section 3. By the
definition of $P_1,P_2$, we can decompose $\Ss_R(P_1)$ into
$\Pi_{>} \Ss_R(P_1)$ and $\sigma_{1} \Ss_R(P_1)$. We can also
decompose $\Ss_R(P_2)^{-1}$ into its restrictions to the images of
$\Pi_{>}$ and $\sigma_{2}$. We denote these maps by
$\Ss_R(P_2)^{-1}\Pi_{>}$ and $\Ss_R(P_2)^{-1}\sigma_{2}$
respectively. Hence the operator
$S_{R,1,2}:=\Ss_R(P_1)\Ss_R(P_2)^{-1}$ has the following form:

$$
\pmatrix{
       \Pi_{>}S_{R,1,2}\Pi_{>} &
\Pi_{>}S_{R,1,2}\sigma_{2} \cr
       \sigma_{1}S_{R,1,2}\Pi_{>} &
\sigma_{1}S_{R,1,2}\sigma_{2}
         \cr}.
$$
\bigskip

\noi Now we see that $\Pi_{>}\Ss_R(P_1)\Ss_R(P_2)^{-1}\Pi_{>}$ is
the identity map on $H^{0}_{>}$ so that it does not depend on $R$.
By definition, $\Pi_{>}\Ss_R(P_1)\Ss_R(P_2)^{-1}\sigma_{2}$ and
$\sigma_{1}\Ss_R(P_1)\Ss_R(P_2)^{-1}\Pi_{>}$ are the zero maps.
Finally we consider the map
$\sigma_{1}\Ss_R(P_1)\Ss_R(P_2)^{-1}\sigma_{2}$. By the Nicolaescu
description of the dynamics of $e^{-R B_{\nu}}L^0_{\nu}$, $e^{-R
B_{0}}L^0_{0}$ does not depend on $R$ so that
$\sigma_{1}\Ss_R(P_1)\Ss_R(P_2)^{-1}\sigma_{2}$ is independent of
$R$. Hence $\Ss_R(P_1)\Ss_R(P_2)^{-1}$ is independent of $R$. Now
the Proposition follows from (\ref{e:sw}).
\end{proof}

\bigskip
\noi
We can combine Proposition \ref{p:R-dep} with the results of \cite{JKPW3}
to obtain
a very interesting result which corresponds to the Lesch--Wojciechowski
formula for
the variation of the $\eta$-invariant.

\bigskip

We have to introduce elements of {\it Scattering Theory} in order
to present the formula. We introduce the manifolds  $M_{2,\infty}$
which are manifolds $M_2$ with the semicylinder $ (-\infty,0]
\times Y$ attached to. Let $\Dd_{2,\infty}$ denote the natural
extension of $\Dd$ to $M_{2,\infty}$. The operator
$\Dd_{2,\infty}$ over $M_{2,\infty}$ has continuous spectrum equal
to $(-\infty,\infty)$. The number $\lambda \in (-\infty,\infty)$
and $\phi\in {\rm ker}(B)$ determine a generalized eigensection of
$\Dd_{2,\infty}$, which has the following form on
$(-\infty,0]\times Y \subset M_{2,\infty}$ (see (4.24) in
\cite{Mu94}):

$$
E(\phi,\lambda)=e^{i\lambda u}(\phi+iG\phi) +e^{-i\lambda
u}C(\lambda)(\phi+iG\phi)
$$
$$
+\theta(\phi,\lambda),
$$

\noi where $\theta(\phi,\lambda)$ is a square-integrable section
of $S$ on $M_{2,\infty}$ which is orthogonal to ${\rm ker}(B)$,
when restricted to $\{u\} \times Y$, and $C(\lambda)$ is the
scattering matrix. We refer to \cite{Mu94} and \cite{Mu98} for the
presentation of the necessary material from {\it Scattering
Theory}.
\bigskip

Let $C : W \to W$ denote a unitary operator acting on the
finite-dimensional vector space $W$. We introduce the operator
$D(C)$ equal to the differential operator $-i\frac12\frac{d}{du}$
acting on $L^2(S^1, E_C)$ where $E_C$ is the flat vector bundle
over $S^1=\mathbb{R}/\mathbb{Z}$ defined by the holonomy
$\overline{C}$.

\bigskip
Now we define the operators
$$
I = (G-i) : {\rm ker} (B) \to {\rm ker} (G+i) \ \ ,
$$
$$
P_{i}=\frac 12(\sigma_i - 1) : {\rm ker}(B) \to {\rm ker}(\sigma_i+1)
$$
and
$$
S_{i}(\lambda)=-P_{i}\circ C(\lambda) \circ I|_{{\rm
ker}(\sigma_i+1)} \ \ .$$

\vskip 5mm

\noi Then $S_{1}:=S_{1}(0)$ and $S_{2}:=S_{2}(0)$ are the unitary
operators acting on the finite-dimensional vector spaces and we
have well-defined self-adjoint, elliptic operators $ D(S_{1}),
D(S_{2})$. The main result of \cite{JKPW3} gives the formula

\begin{equation}\label{e:lim}
\lim_{R\to\infty}
\frac{{\rm det}_{\z}(\Dd_{R})_{P_1}^2}{{\rm det}_{\z}(\Dd_{R})_{P_2}^2}
=\frac{{\rm det}_{\z}D(S_1)^2}{{\rm det}_{\z}D(S_2)^2}
\end{equation}

\bigskip
\noi under the assumption ${\rm ker}(\Dd_R)_{P_1}={\rm ker}(\Dd_R)_{P_2}=0$.
However we showed that the left side of (\ref{e:lim}) is
independent of $R$, hence we have the following Corollary of
Proposition \ref{p:R-dep}:

\bigskip

\begin{cor}\label{c:rel}
Assume that ${\rm ker}(\Dd_R)_{P_1}=$ ${\rm ker}(\Dd_R)_{P_2}=0$.
Then we have
$$
\frac{{\rm det}_{\z}\Dd_{P_1}^2}{{\rm det}_{\z}\Dd_{P_2}^2}
=\frac{{\rm det}_{\z}D(S_1)^2}{{\rm det}_{\z}D(S_2)^2} \ \ .
$$
\end{cor}

\end{document}